\documentclass[12pt, twoside, leqno]{article}

\usepackage{amsfonts,amssymb,amsmath,amscd,amsthm}
\usepackage{graphicx}
\usepackage{enumitem}
\usepackage{hyperref}

\makeatletter

\newtheorem{theorem}{Theorem}[section]
\newtheorem{lemma}[theorem]{Lemma}

\newtheorem{corollary}[theorem]{Corollary}
\newtheorem{observation}[theorem]{Observation}

\theoremstyle{definition}

\newtheorem{conjecture}[theorem]{Conjecture}

\numberwithin{equation}{section}

\newcommand{\pgap}   {{\mathcal G}}
\newcommand {\gap}     {\makebox[0.075 in]{}}   
\newcommand {\biggap}     {\makebox[0.2 in]{}}   
\newcommand {\st}      {\gap : \gap}   
\newcommand{\lil}   {\scriptstyle }

\newcommand {\fto}     {\longrightarrow}

\newcommand {\set}[1]  {\left\{ {#1} \right\}}   
\newcommand {\ord}[1]  {{#1}^{\rm th}}  
\newcommand {\pml}[1]  {{#1}^{\#}}

\begin{document}

\baselineskip=17pt

\title{Discrete dynamics in Eratosthenes sieve}

\author{Fred B. Holt\\
{\tt fbholt@primegaps.info} \\
Seattle, WA}

\date{}

\maketitle

\renewcommand{\thefootnote}{}



\renewcommand{\thefootnote}{\arabic{footnote}}
\setcounter{footnote}{0}

\begin{abstract}
We study Eratosthenes sieve as a discrete dynamic system.  At each stage of the sieve there is a cycle of gaps $\pgap(\pml{p})$ of
length $\phi(p^\#)$ and span $p^\#$.  A recursion $\pgap(p_k^\#)\fto \pgap(p_{k+1}^\#)$ creates the next cycle from
the current one.  

If we take initial conditions from the cycle $\pgap(p_0^\#)$, then for all constellations of span $|s| < 2p_1$, including gaps $g < 2p_1$,
the driving terms of various lengths form Markov chains.  These yield {\it exact} models for the populations $n_s(p_k^\#)$ for all further
stages of the sieve.  

If $s$ is an admissible constellation of length $J$, then its population $n_{s,J}(p^\#)$ grows as $\Theta \left( \prod (q-J-1)\right)$. 
So we factor out the superexponential growth to obtain the
exact model for the relative population $w_s(p_k^\#)$ of the constellation $s$ across all further stages of the sieve.
$$ w_{s,J}(p_k^\#) \; = \; n_{s,J}(p_k^\#) \, / \, \prod_{J+1 < p \le p_k} (p-J-1) $$
The asymptotic value of the relative population is a constant ${w_{s,J}(\infty) \ge 1}$ that depends only on the values $(q-\nu_s(q))$
over the odd prime factors that divide a span in $s$.

Assuming that the instances of a constellation $s$ are approximately uniformly distributed in $\pgap(p_k^\#)$, we develop initial estimates
of the number of instances $s$ that would occur in the {\em interval of survival} $\Delta H(p_k) = (p_k^2, p_{k+1}^2)$.  We define a statistic $\eta_s(p_k)$, 
the quadratic density of the constellation $s$ over the interval $\Delta H(p_k)$.  Our initial
estimates $\widehat{\eta_g}(p)$ for the quadratic densities of prime gaps generally agree with samples up to $5.677\,E14$.
\end{abstract}

\section{Introduction}
We have been studying Eratosthenes sieve as a discrete dynamic system.
At each stage of the sieve, there is a cycle of gaps $\pgap(\pml{p})$
of length $\phi(\pml{p})$ (number of gaps in the cycle) and span $\pml{p}$ (sum of the gaps in the cycle).  

The existence of these cycles is
well-known \cite{Pol} and the early cycles, e.g. $\pgap(\pml{5})$ or $\pgap(\pml{7})$, are routinely rediscovered.  
The cycle $\pgap(p^\#)$ is the cycle of gaps among the $p$-rough numbers, which are the unit $1$ and the remaining candidate primes after
Eratosthenes sieve has advanced through prime $p$.

What is novel about our approach is that we have identified a 3-step recursion that produces the next cycle of gaps from the current one,
$$\pgap(p_k^\#) \fto \pgap(p_{k+1}^\#)$$
The cycles of gaps under this recursion constitute a discrete dynamic system, and we are thereby able to obtain analytic results far beyond what could 
be computed explicitly.

A constellation of length $J$ is a sequence of $J$ consecutive gaps.
Gaps are the special case $J=1$.

\begin{theorem}\label{ThmPops}
Let $s$ be an admissible constellation of length $J$, and let $Q(s)$ be the product of all the odd primes that divide a span between 
boundary fusions in $s$.  Then the aggregate population of $s$ and its driving terms in $\pgap(\pml{p})$ is
\begin{equation}
\sum_{j \ge J} n_{s,j} (\pml{p})  =   \prod_{q \le p} (q - \nu_s(q)). \label{Eqsumns}
\end{equation}
If we have initial conditions $n_{s,j}(p_0^\#)$ in $\pgap(p_0^\#)$ for $p_0$ such that the span $|s| < 2p_1$, 
then the relative population $w_{s,J}(p_k^\#)$ of $s$ in $\pgap(p_k^\#)$ is
\begin{eqnarray}
w_{s,J}(p_k^\#) & = & n_{s,J}(p_k^\#) \gap /   \prod_{J+1 < q \le p_k} (q-J-1) \nonumber \\
	& = & \ell_1 - \ell_2 a_2^k + \ell_3 a_3^k - \cdots +(-1)^{J_1-J} \ell_{J_1-J+1} a_{J_1-J+1}^k \label{Eqwsj}
\end{eqnarray}
with coefficients: \hspace{0.25in} $\ell_j = L_j^T\cdot n_s(p_0^\#) \gap / \gap \prod_{J+1 < q \le p_0} (q-J-1)$ \\[0.08in]
and for parameters the eigenvalues: \hspace{0.2in} $a_j^k = \prod_{p_1\le q \le p_k} \frac{q-J-j}{q-J-1}$.  \\[0.08in]
$J_1$ is the length of the longest driving term for $s$ in $\pgap(p_0^\#)$.

\noindent The asymptotic relative population for $s$, compared to other admissible constellations of length $J$, is
\begin{equation}
w_{s,J}(\infty)  \; = \ell_1 \; =  \prod_{q \le J+1} (q-\nu_s(q))  \cdot  \prod_{\substack{q > J+1, \\ q \mid Q(s)}} \frac{q-\nu_s(q)}{q-J-1}  \label{Eqwsinf}
\end{equation}
\end{theorem}

We develop the proof across Section~\ref{ModelSection} below.

\begin{figure}[hbt]
\centering
\includegraphics[width=5.75in]{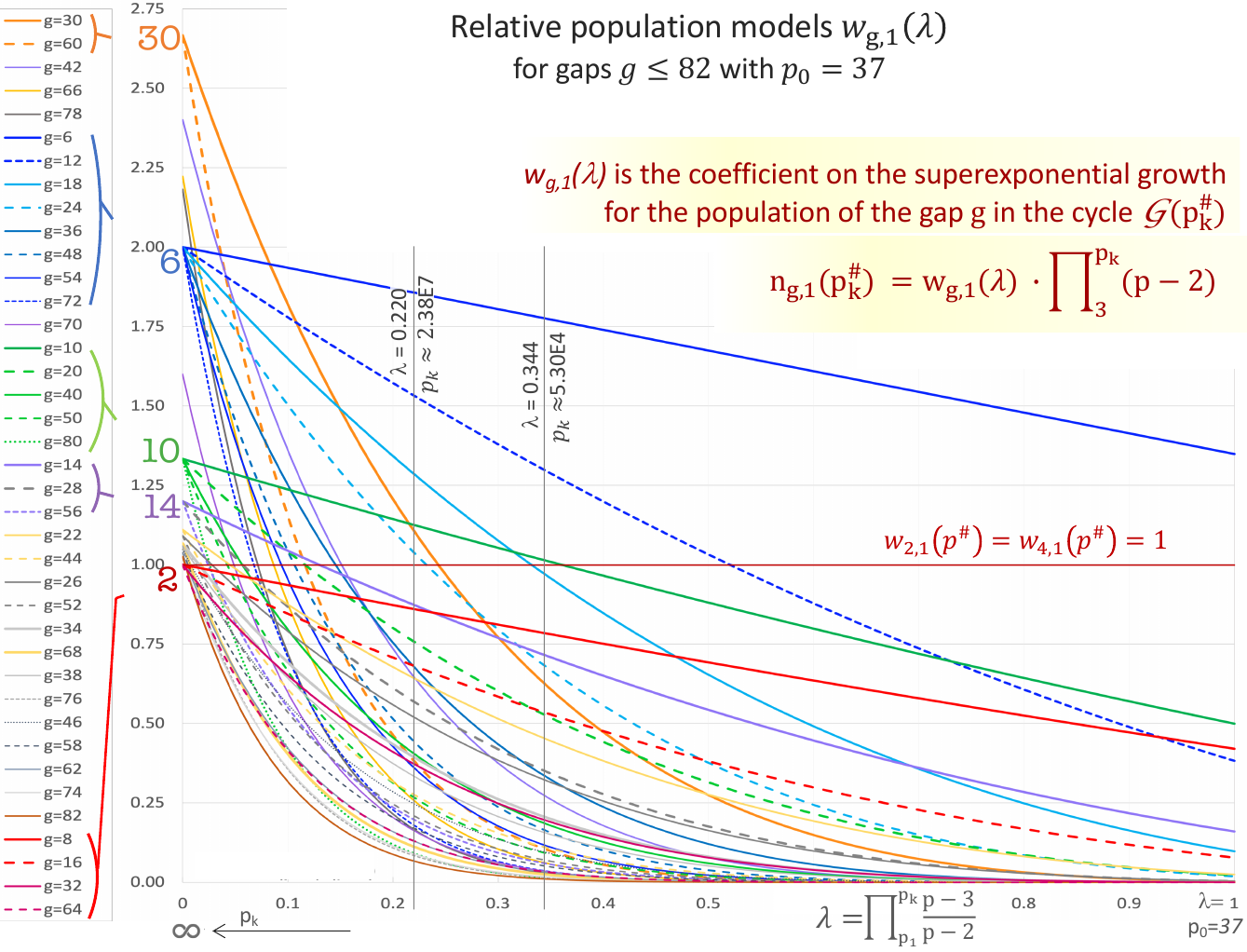}
\caption{\label{AllGapsFig} The relative population models $w_{g,1}(p_k^\#)$ from Equation~\ref{Eqwsj}
for all gaps $2 \le g \le 82$, with initial conditions in $\pgap(37^\#)$.
Each model $w_g(\lambda)$ begins at the right, with $p_0=37$ and the parameter $\lambda=a_2^k=1$.  
As $p_k \fto \infty$ the parameter $\lambda \fto 0$.  The asymptotic relative population of a gap
$g$ depends only on the prime factors of $g$, and we color-code a few families of curves accordingly.}
\end{figure}

These models of the populations $n_{s,J}(p_k^\#)$ and relative populations $w_{s,J}(p_k^\#)$ across the cycles of gaps $\pgap(p_k^\#)$ are exact.
Setting the length $J$, this provides counts of the relative populations of constellations across the cycles $\pgap(p^\#)$ that are consistent with the estimates conjectured
by Hardy and Littlewood \cite{HL}.   For $J=1$ the models are for individual gaps, and the models $w_{g,1}(p_k^\#)$ for gaps $g\le 82$ are shown in Figure~\ref{AllGapsFig}.

Theorem~\ref{ThmPops} tells us that the populations of {\em all} admissible constellations of length $J$ ultimately grow at the same rate, 
$\prod _{p > J+1}(p-J-1)$.  
For each admissible constellation $s$, the asymptotic relative population $w_{s,J}(\infty) \ge 1$ is a constant, which
depends only on the factors $(q-\nu_s(q))$ over the odd prime factors $Q(s)$ of the spans between the boundary fusions in $s$.
So the population $n_{s,J}(p^\#)$ of any admissible constellation of length $J$ grows as 
$$n_{s,J}(p^\#) \; = \; \Theta \left(\prod _{J+1< q \le p}(q-J-1)\right).$$  

\begin{corollary}\label{CorTheta}
Let $s_1, \; s_2$ be admissible constellations of lengths $J_1$ and $J_2$ respectively.
If $J_1=J_2$, the ratio of their populations approaches a positive constant
$$n_{s_2,J}(\pml{p}) \; / \; n_{s_1,J}(\pml{p}) \quad \xrightarrow[p \rightarrow \infty]~\quad w_{s_2,J}(\infty) \; / \; w_{s_1,J}(\infty).$$

\noindent Otherwise wlog let ${J_1 < J_2}$. Then
$$
n_{s_2,J_2}(\pml{p}) \; / \; n_{s_1,J_1}(\pml{p}) \quad \xrightarrow[p \rightarrow \infty]~\quad 0. 
$$
\end{corollary}

Corollary~\ref{CorTheta} holds no matter how populous the longer $s_2$ is, or how rare the shorter $s_1$ is, at any stage of the sieve.

These asymptotic results evolve over incredibly large scales, far beyond the computational horizon.  
We rely on Merten's Third Theorem
\begin{equation}
\prod_{2 \le q \le p} \frac{q-1}{q} \; = \; \frac{e^{-\gamma} + o(1)}{\ln p} \label{EqMerten3}
\end{equation}
to associate large values of $p_k$ with the parameter $\lambda=a_2^k$ for these models $w_{s,J}(p_k^\#)$.
The parameter $\lambda \sim C / \ln p$, for example in the horizontal axis in Figure~\ref{AllGapsFig}.
Analyses of statistical samples such as
\cite{OS} must be considered in light of this evolution \cite{FBHbias}; otherwise the samples can be misinterpreted at face value.

\vspace{0.125in}

The exact models and eigenstructure for the relative populations $w_{s,J}(p_k^\#)$ of small constellations ($|s| < 2p_1$) in the cycles $\pgap(p_k^\#)$
are, we believe, an original contribution.  The population models are consistent with conjectures by Hardy \& Littlewood, but these models $w_s(p_k^\#)$
also capture
the evolution of the populations from longer driving terms to the constellation itself.  
We can calculate the relative population models $w_{s,J}(p_k^\#)$ far beyond cycles $\pgap(p_k^\#)$ or samples of gaps that we can explicitly compute.
The models $w_{s,J}(p_k^\#)$ show that Eratosthenes sieve treats
all admissible constellations of the same length $J$ fairly.  That is, the dynamics do not vary across the admissible constellations of length $J$; 
for a fixed length $J$ the models
vary only in their initial conditions $w_{s,J}(p_0^\#)$.

The proofs for these results are provided below, and accompanying Jupyter notebooks are available at 
\url{https://github.com/fbholt/primegaps-v2}

\section{The dynamic system}\label{SysSection}
At the stage of Eratosthenes sieve in which $p$ has been confirmed as a prime and the multiples of $p$ removed, there is a cycle of gaps
$\pgap(\pml{p})$ among the remaining candidate primes.  These remaining candidate primes are also known as the $p$-rough numbers,
and the first period of the cycle covers the generators of 
$\mathbb{Z} \bmod \pml{p}$.

For example, the early cycles of gaps are
\begin{eqnarray*}
\pgap(\pml{3}) & = & 4 \; 2 \\
\pgap(\pml{5}) & = & 6 \; 4 \; 2 \; 4 \; 2 \; 4 \; 6 \; 2
\end{eqnarray*}
The cycle $\pgap(\pml{p})$ has length $\phi(\pml{p})$ and span $\pml{p}$.

Polignac \cite{Pol} observed these cycles in 1849, as the context for his conjecture about gaps between primes.  
The smaller cycles are used in the wheel factoring method \cite{Wheel}.  Still, it seems that every
few years the cycles $\pgap(\pml{5})$ and $\pgap(\pml{7})$ are rediscovered.

We have not seen any prior exposition about the 3-step recursion \cite{FBHSFU, FBHPatterns} that creates the next cycle of gaps
$\pgap(\pml{p_{k+1}})$ from the current one $\pgap(\pml{p_k})$.

\begin{lemma}\label{R123Lem}
For the cycles of gaps, there is a 3-step recursion that produces $\pgap(\pml{p_{k+1}})$ from 
${\pgap(\pml{p_k}) = g_1 \; g_2 \; \ldots \; g_{\phi(\pml{p_k})}}$.
\begin{itemize}
\setlength{\itemsep}{0.01in}
\item[R1:] {\em Next prime.} $p_{k+1} = g_1 + 1$
\item[R2:] {\em Initial images.} Concatenate $p_{k+1}$ copies of $\pgap(p_k^\#)$.
\item[R3:] {\em Fusions.}  Add together $g_1+g_2$ and thereafter at the running sums indicated by the element-wise product
$p_{k+1} \ast \pgap(\pml{p_k})$.
\end{itemize}
\end{lemma}

\begin{figure}[hbt]
\centering
\includegraphics[width=5.75in]{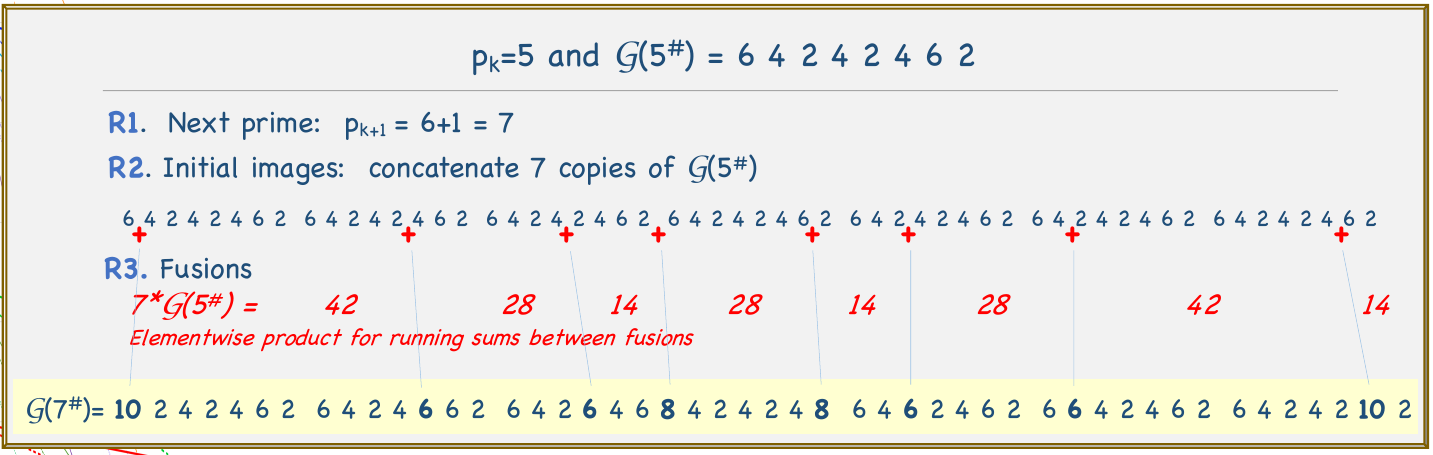}
\caption{\label{R123Fig}  The $3$-step recursion produces the next cycle of gaps $\pgap(p_{k+1}^\#)$ from the current one $\pgap(p_k^\#)$.
We use $\pgap(5^\#)\fto \pgap(7^\#)$ to illustrate this.}
\end{figure}

\begin{proof}
The first gap $g_1$ in $\pgap(p_k^\#)$  is the gap from the unit $1$ to the next-smallest $p$-rough number $\gamma_1$, which is the
next prime $p_{k+1}$.  We concatenate $p_{k+1}$ copies of the current cycle $\pgap(p_k^\#)$, which cover the $p_k$-rough numbers from
$1$ to ${1+ p_{k+1}\cdot p_k^\# = 1+ p_{k+1}^\#}$.  We need only to remove the multiples of $p_{k+1}$ from these initial images to obtain the
$p_{k+1}$-rough numbers covered by the first period of the cycle $\pgap(p_{k+1}^\#)$.

Removing a multiple of $p_{k+1}$ corresponds to fusing the gaps $g_i+g_{i+1}$ on either side. The $p_k$-rough numbers that are not $p_{k+1}$-rough
numbers have the form $p_{k+1} \cdot \gamma$ where $\gamma$ is a $p_k$-rough number.  The difference between multiples of $p_{k+1}$
with consecutive $p_k$-rough numbers 
$$ d \, = \, p_{k+1}\cdot \gamma_{i+1} - p_{k+1} \cdot \gamma_{i} = p_{k+1} \cdot (\gamma_{i+1}-\gamma_i) \, = \, p_{k+1} \cdot g_{i+1}$$
Thus the elementwise product $p_{k+1} \ast \pgap(p_k^\#)$ is the sequence of differences from one remaining multiple of $p_{k+1}$ to the next.

We note that the cycle $\pgap(p_k^\#)$ always ends with a gap $g=2$, and this elementwise product $2 p_{k+1}$ wraps around the end of the
initial images, back to the first fusion $g_1+g_2$.
\end{proof}

\begin{figure}[hbt]
\centering
\includegraphics[width=4.5in]{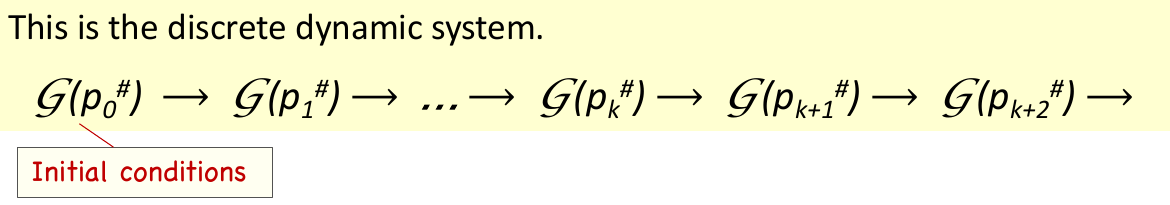}
\caption{\label{DynSysFig}  The cycles of gaps $\pgap(p^\#)$ under the 3-step recursion form a discrete dynamic system. If we take initial populations 
$n_{s,j}(p_0^\#)$ for constellations and their driving terms from $\pgap(p_0^\#)$, then we can form exact models $w_{s}(p_k^\#)$ for the relative
population of every 
constellation $s$ with span $|s|\le 2p_1$.}
\end{figure}

The cycles of gaps $\pgap(\pml{p})$ under this recursion form a discrete dynamic system.  We do not need any external information.
Each cycle of gaps $\pgap(p^\#)$ contains all the information we need to create all subsequent cycles.

The value in this dynamic system is not the construction of complete cycles $\pgap(p_k^\#)$.  There are more efficient methods for finding
$p$-rough numbers, and these cycles have lengths $\phi(p_k^\#)$ which very quickly exhaust computational resources.  
The valuable insight is that the dynamic system preserves a lot of structure from ${\pgap(p_k^\#)\rightarrow \pgap(p_{k+1}^\#)}$, and constellations
evolve in predictable ways \cite{FBHSFU, FBHPatterns, FBHktuple}.

Two important observations facilitate creating the population models for constellations of gaps.

\begin{observation}\label{Obs1}
In step R3 of the recursion $\pgap(p_k^\#) \fto \pgap(p_{k+1}^\#)$ each of the $\phi(p_k^{\#})$ possible fusions 
in $\pgap(p_k^\#)$ occurs exactly once.
\end{observation}

This is a result of the Chinese Remainder Theorem in this context.  Lemma~\ref{LemAdmIter} provides details supporting this observation.

\begin{observation}\label{Obs2}
In step R3 of the recursion the minimum distance between fusions is $2 p_{k+1}$.
\end{observation}

\begin{figure}[hbt]
\centering
\includegraphics[width=3.1in]{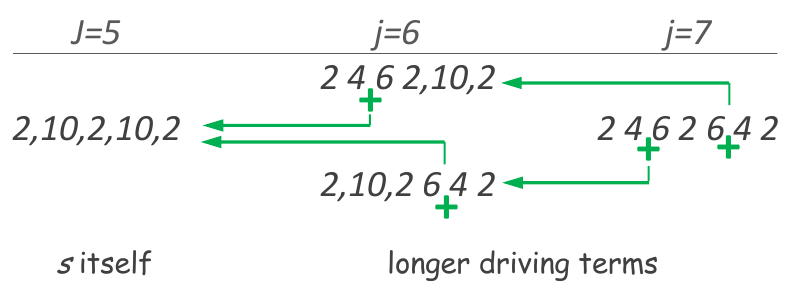}
\caption{\label{J5DynFig} For $s=2,10,2,10,2$ of length $J=5$, we show its admissible driving terms, and we mark its interior fusions
under the recursion.}
\end{figure}

There may be longer constellations $\tilde{s}$ that form $s$ under the fusions of step R3 in the recursion.  For example, if we consider
$s=2,10,2,10,2$ of length $J=5$, the constellation $\tilde{s}=2462642$ produces $s$ when the pairs of adjacent gaps $46$ and $64$ 
are fused into $10$'s.
We call these constellations {\em driving terms} for $s$.  This is illustrated in Figure~\ref{J5DynFig}.

For a driving term $\tilde{s}$ of length $j \ge J$ for $s$, the $J+1$ fusions that would eliminate the image from being a driving term for 
$s$ are the {\em boundary fusions} for $\tilde{s}$.  The $j-J$ fusions for which the image of $\tilde{s}$ is still a driving term
for $s$ are the {\em interior fusions} for $\tilde{s}$.  For $s$ itself all $J+1$ of the fusions are boundary fusions.

\begin{lemma}
Let $s$ be an admissible constellation of length $J$ in the cycle $\pgap(p_0^\#)$.  If the span $|s| < 2 p_{1}$, then all $J+1$ possible fusions
in $s$ in step R3 occur in different initial images of $s$ in step R2.
\end{lemma}

{\bf Proof.} This is a direct result of the two observations.
$\square$

This condition $|s| < 2p_1$ suffices to put the relative populations of $s$ and its driving terms in $\pgap(p_0^\#)$ into a Markov chain, and the associated
linear system has a beautifully simple eigenstructure  \cite{FBHSFU, FBHPatterns}, from which we can exactly model the populations in all
subsequent cycles of gaps $\pgap(p_k^\#)$.  We develop the models for the relative populations $w_s(p^\#)$ in Section~\ref{ModelSection} below.

Before moving onto those exact population models, let's look more closely at steps R2 and R3 of the recursion and their ties to admissibility.

\subsection{Admissible constellations}
A $k$-tuple is given in terms of generators, and a constellation is given in terms of the gaps between successive generators.  
A constellation of length $J$ corresponds to a family of $k$-tuples with ${k = J+1}$.  

An {\em instance} of a constellation
$s$ is identified by specifying the initial generator $\gamma_0$.  Once $\gamma_0$ is set, the other generators in the $k$-tuple are 
determined by ${\gamma_j = \gamma_{j-1} + g_j}$.

At this point we focus on boundary
fusions for an instance of $s$ or one of its driving terms.  
We define a {\em general instance} of $s$ to be an instance of $s$ itself or an instance of a longer driving term for $s$.
A boundary fusion eliminates a general instance from being a driving term for $s$.  

For a constellation $s$ of length $J$ the parameter $\nu_s(p)$ is the number of residue classes $\bmod~p$ covered by 
the boundary fusions of any general instance
of $s$ \cite{HL}, with 
$$1 \le \nu_s(p) \le \min\set{p, J+1}.$$

If $\tilde{s}$ is a longer driving term for $s$ of length $j > J$, then $\tilde{s}$ has $J+1$ boundary fusions in this role.
The remaining $j-J$ fusions in $\tilde{s}$ are {\em interior fusions}.  An interior fusion shortens the image of $\tilde{s}$ but does not
change the image's status as a driving term for $s$.

We call $r=\gamma_0 \bmod p$ an {\em admissible instance} of $s$ in $\pgap(\pml{p})$ iff
${\gamma_j \neq 0 \bmod p}$ for all ${0 \le j \le J}$.
We use $s$ as the reference for admissibility and boundary fusions.
An admissible instance of $s$ is a general instance of $s$, for which
$\gamma_0$ preserves the boundary fusions.

There are $p-\nu_s(p)$ admissible instances $\gamma_0$ for $s$, modulo $p$. 
Define $\Upsilon$ to be the set of admissible residues for the initial generator $ \gamma \bmod p$.
$$\Upsilon_s(p) = \set{ r = \gamma \bmod p \st \gamma_0=\gamma \; {\rm is~admissible~for} \; s}$$
The cardinality of $\Upsilon_s(p)$ is $p-\nu_s(p)$.

\begin{table}[h]
\centering
\begin{tabular}{rccc} \hline
\multicolumn{4}{c}{Admissibility for $s=2,10,2,10,2$.} \\
\multicolumn{1}{c}{$p$} & $\nu_s(p)$ & $\Upsilon_s(p) $ & $p-\nu_s(p)$\\ \hline
$5$ & $4$ & $\set{2}$ & $1$\\
$7$ & $4$ & $\set{1,3,6}$ & $3$ \\
$11$ & $5$ & $\set{1,2,3,4,5,6}$ & $6$ \\
$13$ & $5$ & $\set{3,4,5,6,7,8,9,10}$ & $8$\\
$17$ & $6$ & $\set{1,2,4,6,7,9,11,12,13,14,16}$ & $11$ \\
\end{tabular}
\caption{Illustrating $\Upsilon_s(p)$ for $s=2,10,2,10,2$.  Once $p_k \ge 13$, $|s|<2p_{k+1}$ so all fusions occur in separate images under step R2, and ${\nu_s(p_{k+1})=6}$.}
\end{table}

We say that a constellation $s$ is {\em admissible for $p$} iff $\nu_s(p) < p$, or equivalently iff $\Upsilon_s(p)$ is non-empty.  
By default the constellation $s$ is admissible for all primes $p > J+1$.

Finally, a constellation $s$ is {\em admissible} iff $s$ is admissible for all primes $p$.

\vspace{0.1in}

Let $Q(s)$ denote the product of odd primes that divide a span between boundary fusions in $s$.  For the example $s=2,10,2,10,2$ this product
${Q(s)=3\cdot 5 \cdot 7 \cdot 11 \cdot 13}$.

\vspace{0.125in}

\noindent
\begin{minipage}{\linewidth}
\begin{lemma}\label{Lemnu}
Let $s$ be an admissible constellation of length $J$, and let $p$ be an odd prime.
\begin{enumerate}
\setlength{\itemsep}{0.01in}
\item[a.] If $p \le J+1$, then $p \mid Q(s)$. 
\item[b.] If $p > J+1$, then $\nu_s(p) < J+1$ iff $p \mid Q(s)$.
\end{enumerate}
\end{lemma}
\end{minipage}

\vspace{0.125in}

{\bf Proof.} 
This is a remark on the top of p.62 in \cite{HL}.
For (a), when $p \le J+1$ there are at least as many boundary fusions as residue classes. Since $s$ is admissible, there must be multiple fusions
in some residue classes, and this requires ${p \mid Q(s)}$.

When $p > J+1$ there are more residue classes than boundary fusions, and the admissibility of $s$ is automatic.  
In this case, we note that 
$$ p - \nu_s(p) > p-J-1 \gap {\rm iff} \gap p \mid Q(s).$$
If all $J+1$ boundary fusions occur in separate images of $s$, then ${\nu_s(p)=J+1}$.  Otherwise some boundary fusions occur in the
same image of $s$.  This requires ${p \mid Q(s)}$, and $\nu_s(p) < J+1$.
$\square$

\subsection{Admissible instances across Eratosthenes sieve}
Under the recursion, Eratosthenes sieve systematically produces all admissible instances of every admissible constellation.
The admissible instances of $s$ include its driving terms.

We begin with a cycle of gaps $\pgap(p_0^\#)$ for a $p_0$ that is manageable.  If needed, we could begin with the trivial starting points
$p_0=3$ or $5$, but we often start with ${7 \le p_0 \le 37}$.

\begin{figure}[hbt]
\centering
\includegraphics[width=4.5in]{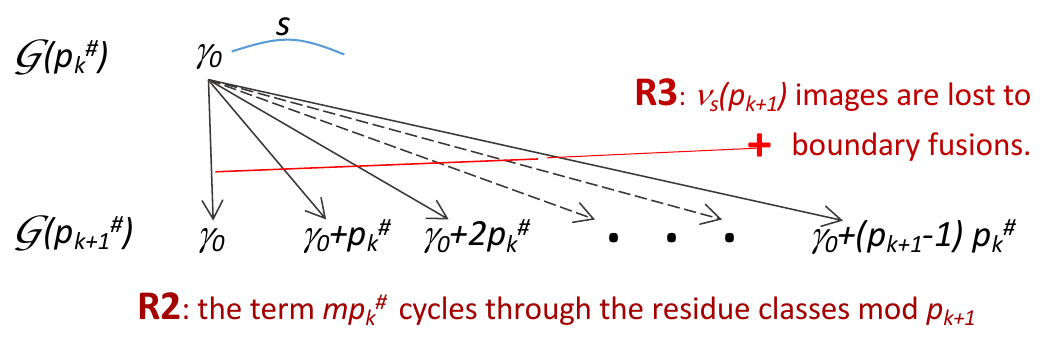}
\caption{\label{RFig} Under the recursion each instance $\gamma_0$ of $s$ in $\pgap(\pml{p_k})$ is replicated $p_{k+1}$ times
at offsets of $\pml{p_k}$.  Then $\nu_s(p_{k+1})$ of these images are lost to boundary fusions and $p_{k+1}-\nu_s(p_{k+1})$ survive
in $\pgap(p_{k+1}^\#)$. }
\end{figure}

\begin{lemma}\label{LemAdmIter}
Let $\gamma_0$ be the initial generator for an admissible instance of $s$, possibly a longer driving term for $s$, in $\pgap(p_0^{\#})$.
Let
\begin{equation}
\gamma_{0,k} = \gamma_0 \; + \; m_1 \cdot p_0^{\#} \; + \; m_2 \cdot p_1^{\#}  + \cdots + m_k \cdot p_{k-1}^{\#} \label{EqGamk}
\end{equation}
be an admissible instance of $s$ in $\pgap(\pml{p_k})$, with $0 \le m_j < p_j$.  
Then in step R2 of the recursion
${\pgap(\pml{p_k}) \fto \pgap(\pml{p_{k+1}})}$, this general instance $\gamma_{0,k}$ of $s$ has $p_{k+1}$ initial images, with initial generators
$$ \gamma_{0,k} \; + m \cdot p_k^\# \gap = \gap \left(\gamma_0 \; + m_1 \cdot \pml{p_0} \; +  \cdots + m_k \cdot p_{k-1}^{\#}\right) + m \cdot p_k^{\#} $$
for $m= 0,1,\ldots,p_{k+1}-1$.  These initial images cycle through all of the residue classes $\bmod \; p_{k+1}$.
As such, ${p_{k+1}-\nu_s(p_{k+1})}$ of these survive in $\pgap(p_{k+1}^{\#})$, and $\nu_s(p_{k+1})$ are eliminated by boundary fusions.
\end{lemma}

{\bf Proof.}
The replication $\gamma_{0,k}+m\cdot p_k^{\#}$ of this general instance of $s$ is induced by step R2 of the recursion.
We have to show that these images $\gamma_{0,k}+m \cdot p_k^{\#}$ cycle through the residue classes $\bmod~p_{k+1}$,
and that $p_{k+1}-\nu_{p_{k+1}}(s)$ of these survive the fusions in step R3.

Since $\pml{p_k} \neq 0 \bmod~p_{k+1}$, for $0 \le m < p_{k+1}$ the initial generators ${\gamma_{0,k} + m \cdot p_{k}^{\#}}$ cycle once through the residue classes
modulo $p_{k+1}$.

For each $m$, $\gamma_{0,k}+m\cdot p_k^{\#}$ survives the fusions in step R3 iff
$$(\gamma_{0,k}+ m\cdot p_k^{\#} ) \bmod p_{k+1} \; \in \; \Upsilon_s(p_{k+1}).$$  
These fusions correspond to removing
multiples of $p_{k+1}$ from $\pgap(p_{k+1}^{\#})$, and each such multiple is ${0 \bmod p_{k+1}}$.
The cardinality of $\Upsilon_s(p_{k+1})$ is ${p_{k+1}-\nu_s(p_{k+1})}$, and we have our result.
$\square$

\begin{figure}[hbt]
\centering
\includegraphics[width=3.25in]{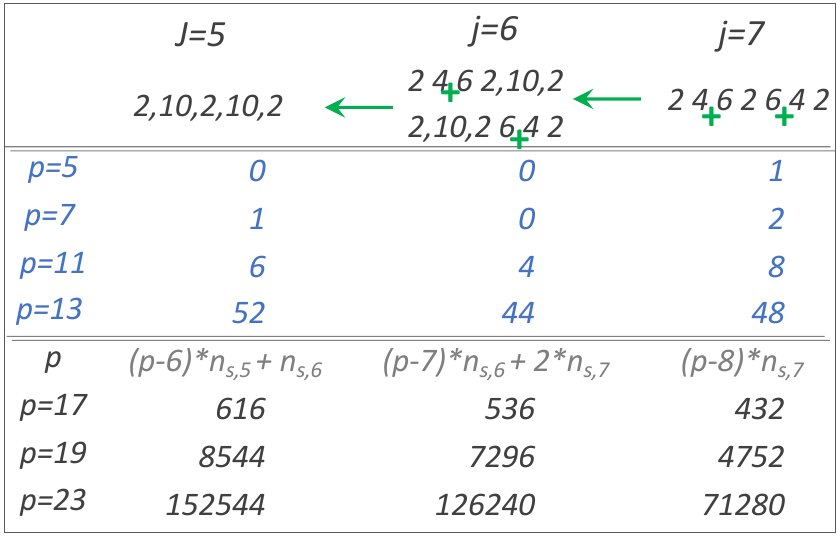}
\caption{\label{segFig} The populations $n_{s,j}$ of $s=2,10,2,10,2$ and its driving terms in $\pgap(p^\#)$ across the first few primes.
From $p=13$ on, the condition ${|s|<2p_1}$ is satisfied and the linear system holds. }
\end{figure}

Eratosthenes sieve methodically eliminates all inadmissible instances of a constellation $s$, and the sieve methodically produces all admissible 
instances of $s$.  As $m$ advances across $0 \le m < p_{k+1}$, every constellation advances through the residues $\bmod~p_{k+1}$ in the same order, with
the stride set by $\Delta r = p_k^\# \bmod p_{k+1}$; two constellations in $\pgap(p_k^\#)$ may differ in their starting points $\gamma_{0,k} \bmod p_{k+1}$, 
but their initial images travel through the residues modulo $p_{k+1}$ in the same order.
If a constellation $s$ is not admissible for the prime $p_{k+1}$, then every general instance of $s$ is eliminated by
the recursion from ${\pgap(p_k^\#) \fto \pgap(p_{k+1}^\#)}$.

Our focus in this section is on the total population of general instances of $s$ that occur in $\pgap(p_{k+1}^\#)$.  Under the recursion, every instance of
a longer driving term for $s$ in $\pgap(p_k^\#)$ is either (a) eliminated by a boundary fusion, (b) replicated as it is (no fusions), or (c) shortened
as a driving term for $s$ in $\pgap(p_{k+1}^\#)$.  For each tuple $(\gamma_0, r_1, r_2, \ldots, r_k)$ of admissible instances $\gamma_0$ in 
$\pgap(p_0^\#)$ and admissible residues $r_i \bmod p_i$ for $s$, there is a unique
general instance of $s$ in $\pgap(p_k^\#)$.

Lemma~\ref{LemAdmIter} shows that each admissible instance of $s$ 
in $\pgap(p_k^\#)$ has ${p_{k+1}-\nu_s(p_{k+1})}$ admissible images in $\pgap(p_{k+1}^\#)$.
Formalizing this, a further lemma provides the population count for admissible instances of a constellation across stages
of Eratosthenes sieve, Equation~\ref{Eqsumns} of Theorem~\ref{ThmPops}.

\begin{lemma}\label{LemNCnt}
Let $s$ be a constellation of length $J$ admissible through prime $p$.  Let $n_{s,j}(p^\#)$ denote the population of admissible driving terms for $s$ of length
$j$ in the cycle of gaps $\pgap(p^\#)$.  Then
\begin{equation}
\tag{\ref{Eqsumns}}
\sum_{j \ge J} n_{s,j}(p^\#) \gap = \gap \prod_{q \le p} (q - \nu_s(q)) 
\end{equation}
\end{lemma}

{\bf Proof.} This follows directly from the previous Lemma~\ref{LemAdmIter}. 
At each stage of the sieve every admissible instance of $s$ has the same number of admissible images under the recursion.  The lefthand side of
Equation~\ref{Eqsumns} adds up the populations of driving terms of $s$ across all lengths $j \ge J$, and the righthand side applies the result from
Lemma~\ref{LemAdmIter}.
$\square$

\section{Exact models of relative populations $w_{s,J}(\lambda)$}\label{ModelSection}
\begin{figure}[hbt]
\centering
\includegraphics[width=5.5in]{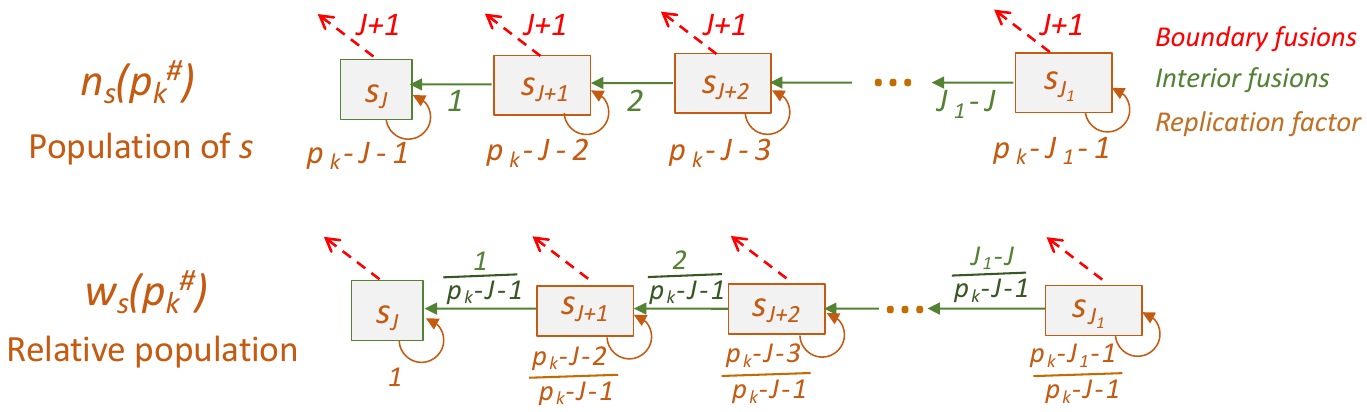}
\caption{\label{DynFig} Under the recursion, once $|s| < 2 p_1$ then the relative populations of $s$ and its driving terms of length $j$ form
a Markov chain for all subsequent $p_k$.  In creating the cycle $\pgap(p_k^\#)$ the dynamics for all admissible constellations of length $J$ 
and span  $|s|< 2 p_k$ are identical.}
\end{figure}

Under the dynamics of the Markov chain depicted in Figure~\ref{DynFig}, 
the populations shift from the longer driving terms to copies of $s$ itself.  When $|s| < 2p_1$ we are able to calculate the exact 
relative population model $w_{s,J}(p_k^\#)$ for all $p_k \ge p_0$.

We denote the population of the constellation $s$ in the cycle of gaps $\pgap(\pml{p})$ as $n_{s,J}(\pml{p})$, and $n_{s,j}(\pml{p})$
denotes the population of the driving terms $\tilde{s}$ of length $j \ge J$.

These populations $n_{s,J}$ grow super-exponentially, dominated by factors of ${(q-J-1)}$.
We define the {\em relative population} $w_{s,J}(\pml{p})$  by
\begin{equation}
w_{s,J}(\pml{p}) = n_{s,J}(\pml{p}) \gap / \prod_{J+1 < q \le p} (q-J-1). \label{Eqwn}
\end{equation}
This Equation~(\ref{Eqwn}) isolates the dominant factors ${(q-J-1)}$.  If we take initial conditions in a cycle $\pgap(p_0^\#)$
for which ${|s| < 2p_1}$, then the relative population $w_{s,J}(\pml{p})$ 
is an iterative linear system as depicted in Figure~\ref{DynFig}.

\begin{figure}[hbt]
\centering
\includegraphics[width=6in]{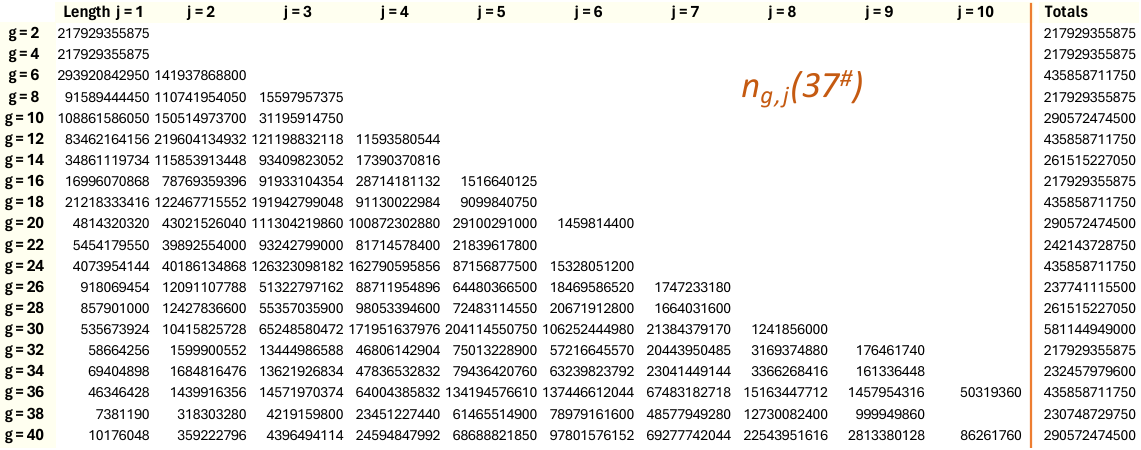}
\caption{\label{G37Fig} A sample of the initial conditions $n_{g,1}(p_0^\#)$ with $p_0=37$.}
\end{figure}

The transfer matrix $M_J(\pml{p})$ at each stage is a banded matrix with diagonal elements $\frac{p-J-j}{p-J-1}$ and
super-diagonal elements $\frac{j}{p-J-1}$.  We then have
$$
w_s(p_{k}^\#)  = M_J(p_k) \cdot w_s(p_{k-1}^\#) \;  = \; M_J^k \cdot w_s(p_{0}^\#),
$$
using the notation
\begin{eqnarray*}
M_J^k & = & M_J(p_k)  \cdot M_J(p_{k-1}) \cdots  M_J(p_1) \\
    & = & R \cdot \Lambda_J^k \cdot L^T.
 \end{eqnarray*}
The matrix of right eigenvectors $R$ is an upper triangular Pascal matrix of alternating sign.  The matrix of left eigenvectors
$L^T$ is an upper triangular Pascal matrix.

The matrix of eigenvalues $\Lambda^k$ is
\begin{eqnarray*}
\Lambda_J^k & = & {\rm diag} \left(1, \gap \prod_{p_1}^{p_k}\frac{p-J-2}{p-J-1}, \gap \prod_{p_1}^{p_k}\frac{p-J-3}{p-J-1}, \gap \ldots, 
  \gap \prod_{p_1}^{p_k} \frac{p-J_1-1}{p-J-1} \right) \\
 & = & {\rm diag}\left( 1,  \gap \lambda=a_2^k, \gap a_3^k, \ldots , \gap a_{J_1-J+1}^k \right)
\end{eqnarray*}
where $J_1 \ge J$ is the length of the longest admissible driving term for $s$.

To verify the eigenstructure for $M_J(p)$, we can proceed by induction on the size of the matrix.  
For the base case use $n_0=3$ or $n_0=4$. When $M_J$ is $n \times n$
the last two entries in the $\ord{j}$ row of $L^T$ are $\binom{n-2}{j-1}$ and $\binom{n-1}{j-1}$.  To declutter the calculations, let 
$\beta = \frac{1}{p-J-1}$. The last column in $M_J$ consists
of $(n-2)$ zeroes, then an entry of $(n-1)\beta$ and an entry of $1-(n-1)\beta$. To show that $L^T\cdot M_J = \Lambda_J \cdot L^T$,
we assume that this holds for the first $(n-1)$ columns of the product.  We need only confirm that for each row $L_j^T$, its product with the
last column of $M_J$ is ${(1-(j-1)\beta)\binom{n-1}{j-1}}$ as $1\le j \le n$.  For $1\le j < n$ we have
\begin{eqnarray*}
L_j^T \cdot \left[ \begin{array}{c} \lil \langle 0 \rangle \\ \lil (n-1)\beta \\ \lil 1-(n-1)\beta \end{array} \right] & = & \binom{n-2}{j-1} (n-1)\beta + \binom{n-1}{j-1}(1-(n-1)\beta) \\
 & = & \binom{n-1}{j-1} (n-j)\beta  + \binom{n-1}{j-1} (1-(n-1)\beta) \\
 & = & (1-(j-1)\beta)\binom{n-1}{j-1}
\end{eqnarray*}
For $j=n$, the $\ord{n}$ row has a single nonzero entry of $1$ on the diagonal,
and the product is $1\cdot (1-(n-1)\beta)$.

From the top row of
\begin{equation}
 w_s(p_k^\#) = M_J^k \cdot w_s(p_0^\#) \; = \; R \cdot \Lambda^k \cdot L^T \cdot w_s(p_0^\#), \label{Eqw}
\end{equation}
we extract an equation for the relative population of $s$ itself:
\begin{equation}
\tag{\ref{Eqwsj}}
w_{s,J}(p_k^\#) = \ell_1 \; - \; \ell_2 a_2^k \; + \; \ell_3 a_3^k - \cdots + (-1)^{j+1} \ell_{j} a_{j}^k + \cdots
\end{equation}
with coefficients $\ell_j = L_j^T\cdot w_s(p_0^\#)$.  This is Equation~\ref{Eqwsj} in Theorem~\ref{ThmPops}.

\subsubsection{The parameter $\lambda=a_2^k$.} 
The dominant eigenvalue in $\Lambda_J^k$ is $1$.  The other eigenvalues
$$
a_j^k \; = \; \prod_{p_1}^{p_k} \frac{p-J-j}{p-J-1}
$$
are all positive, and all converge to $0$ as $p_k \fto \infty$.  We set $\lambda = a_2^k$, the subdominant eigenvalue, as our system parameter.  

The other eigenvalues $a_j^k$ are dominated by $a_2^k$. For $k \ge 1$
$$ 1 \; > \; a_2^k \; > \;  a_3^k \; > \; a_4^k \cdots > 0.$$
In fact $a_j^k \asymp \lambda^{j-1}$ with $0 < a_j^k \le \lambda^{j-1}$.
Thus we can view the evolution of the model $w_{s,J}(p_k^\#)$ through the
evolution of the parameter $\lambda=a_2^k$.  See Figure~\ref{AllGapsFig}.

\begin{table}[h]
\centering
\begin{tabular}{r|rrrrrr} \hline
\multicolumn{7}{c}{Ratios $a_j^k / \lambda^{j-1}$ when $p_0 = 37$} \\
\multicolumn{1}{c}{$J$} & \multicolumn{1}{c}{$a_3^k/ \lambda^2$} & \multicolumn{1}{c}{$a_4^k/ \lambda^3$} & \multicolumn{1}{c}{$a_5^k/ \lambda^4$} &
\multicolumn{1}{c}{$a_6^k/ \lambda^5$} & \multicolumn{1}{c}{$a_7^k/ \lambda^6$} & \multicolumn{1}{c}{$a_8^k/ \lambda^7$}  \\ \hline
1 & \footnotesize 0.994043 & \footnotesize 0.982049 & \footnotesize 0.964043 & \footnotesize 0.940159 & \footnotesize 0.910644 & \footnotesize 0.875855 \\
2 & \footnotesize 0.993854 & \footnotesize 0.981478 & \footnotesize 0.962898 & \footnotesize 0.938258 & \footnotesize 0.907822 & \footnotesize 0.871971 \\
3 & \footnotesize 0.993654 & \footnotesize 0.980871 & \footnotesize 0.961682 & \footnotesize 0.936240 & \footnotesize 0.904827 & \footnotesize 0.867850 \\
4 & \footnotesize 0.993440 & \footnotesize 0.980226 & \footnotesize 0.960388 & \footnotesize 0.934093 & \footnotesize 0.901641 & \footnotesize 0.863470 \\
5 & \footnotesize 0.993213 & \footnotesize 0.979538 & \footnotesize 0.959009 & \footnotesize 0.931804 & \footnotesize 0.898247 & \footnotesize 0.858807 \\
6 & \footnotesize 0.992970 & \footnotesize 0.978803 & \footnotesize 0.957535 & \footnotesize 0.929359 & \footnotesize 0.894624 & \footnotesize 0.853833 \\
7 & \footnotesize 0.992710 & \footnotesize 0.978016 & \footnotesize 0.955958 & \footnotesize 0.926743 & \footnotesize 0.890749 & \footnotesize 0.848516 \\
8 & \footnotesize 0.992432 & \footnotesize 0.977171 & \footnotesize 0.954265 & \footnotesize 0.923936 & \footnotesize 0.886593 & \footnotesize 0.842821 \\
\end{tabular}
\caption{The ratios $a_j^k / \lambda^{j-1}$ for $3 \le j \le 8$ and constellation lengths $1 \le J \le 8$.  Here $p_0 = 37$ and $p_k=547709$.}
\end{table}

In Figure~\ref{AllGapsFig} we plot the relative population models $w_{g,1}(\lambda)$ for the gaps $2 \le g \le 82$.  Here $J=1$
and $p_0=37$.  Taking initial conditions $n_{g,j}(37^\#)$ --- a sample is shown in Figure~\ref{G37Fig} --- we can calculate the coefficients
$\ell_j$ for each gap $g < 82$.
\begin{eqnarray*}
w_{g,j}(p_0^\#) & = & n_{g,j}(p_0^\#) / \prod_3^{p_0} (p-2) \\
\ell_j & = & L^T_j \cdot w_g(p_0^\#).
\end{eqnarray*}
These are the coefficients for Equation~\ref{Eqwsj} for the gaps $g < 82$.

The relative populations for the gaps evolve from right to left with parameter $\lambda \; = \; a_2^k \; = \; \prod_{p_1}^{p_k} \frac{q-3}{q-2}$.

For $p_0=37$ the parameter $\lambda = 1$. As $p_k$ grows, 
Merten's Third Theorem tells us that 
$$\lambda = \Theta \left(\frac{1}{\ln p_k}\right) \; \xrightarrow[p \rightarrow \infty] ~\gap 0.$$
From this relation we see that halving $\lambda$ corresponds to primes approximately the square of $p_k$.
\hspace{0.05in} If $\lambda = \lambda(p_k)$, then $\lambda(q) \approx \frac{\lambda}{2}$ iff $q \approx p_k^2$.

An important observation here is that the system dynamics, e.g. the eigenvectors and eigenvalues, depend only on the primes and the length $J$.  
They do not vary by the specific
constellation of length $J$.  The various admissible constellations of length $J$ differ only in their initial conditions $w_s(p_0^\#)$, including the
length of their longest admissible driving term.

\subsubsection{The case $|s| = 2p_1$.}
We are also able \cite{FBH2p} to develop the model $w_{s,J}(p_k^\#)$ for the case $|s|= 2p_1$, but the first iteration of this case 
requires special handling.  

Let $s$ be an admissible constellation of span $|s|=2p_1$ in $\pgap(p_0^\#)$.  Then under the first recursion ${\pgap(p_0^\#)\fto \pgap(p_1^\#)}$,
for every driving term of $s$ in $\pgap(p_0^\#)$ its first and last boundary fusions occur in a single initial image, and all other fusions occur in separate
initial images.  Observation~\ref{Obs2} guarantees this.

\begin{figure}[hbt]
\centering
\includegraphics[width=4.75in]{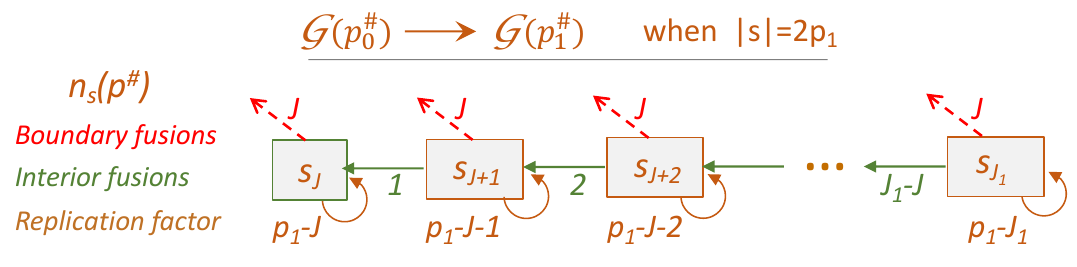}
\caption{\label{2p1Fig} If $|s|=2p_1$, the first iteration has the transition weights as shown. These differ from the transition weights in
Figure~\ref{DynFig} since the first and last boundary fusions occur in a single initial image during the recursion ${\pgap(p_0^\#)\fto \pgap(p_1^\#)}$. }
\end{figure}

The parameter $\lambda$ depends on the initial prime $p_0$.  Thus if we want to keep the models $w_{s,J}(p_k^\#)$ aligned on $p_0$ and the associated
$\lambda$, we can advance the model for the constellation $s$ one step, into $\pgap(p_1^\#)$ and then pull the population vector $n_{s}(p_1^\#)$ back into 
$\pgap(p_0^\#)$ by applying the inverse matrix $(M_J^1)^{-1}$.

This gives us an equivalent surrogate starting point in $\pgap(p_0^\#)$ for $s$ when $|s|=2p_1$.  The true initial population vector in $\pgap(p_0^\#)$ will disagree with
the surrogate initial population vector, but the populations will agree from $\pgap(p_1^\#)$ onward.

We include the graph for $g=82$ in Figure~\ref{AllGapsFig}, but technically this gap does not fall under the conditions for Theorem~\ref{ThmPops}.
The model $w_{82,1}(\lambda)$ is generated using the method described here:  finding initial conditions $n_{82}(37^\#)$, calculating $n_{82}(41^\#)$, then calculating
a surrogate starting point $(M_1^1)^{-1} \cdot n_{82}(41^\#)$, and using the surrogate starting point to calculate the coefficients $\ell_j$ for the model.

Given a cycle $\pgap(p_0^\#)$ from which we can extract initial conditions $n_{s}(p_0^\#)$, we can derive exact models $w_{s}(p_k^\#)$ for all constellations
with $|s| \le 2p_1$.  For $|s| < 2p_1$ we use Theorem~\ref{ThmPops} directly, and for $|s| = 2p_1$ we identify a surrogate starting point from which to derive
the coefficients $\ell_j$.  The models are all exact from $p_1$ onward.  The model for $|s|=2p_1$ will offer the surrogate starting point in $\pgap(p_0^\#)$
rather than the actual populations $n_s(p_0^\#)$, but we rarely use the models to obtain the initial conditions.

\subsection{Asymptotic values}

As $p_k \fto \infty$, the parameter $\lambda \fto 0$, and the relative population model~\ref{Eqwsj} converges to
${w_{s,J}(\infty) = \ell_1 = L_1^T\cdot w_s(p_0^\#)}$.  This holds for any prime
$p_0$ such that $|s| < 2p_1$.
The first left eigenvector is all $1$'s, so 
\begin{eqnarray*} 
w_{s,J}(\infty) =  \ell_1 & = & L_1^T \cdot w_s(p_0^\#) = \sum_J^{J_1} w_{s,j}(p_0^\#)  \\
  & = & \left( \sum_J^{J_1} n_{s,j}(p_0^\#) \right) / \prod_{J+1 < q \le p_0} (q-J-1)  \\
  & = &  \prod_{q \le J+1 } (q-\nu_q(s)) \cdot  \prod_{J+1 < q \le p_0} \frac{q-\nu_q(s)}{q-J-1}   \\
  & = &  \prod_{q \le J+1 } (q-\nu_q(s)) \cdot  \prod_{\substack{J+1 < q \\  q | Q(s)}} \frac{q-\nu_q(s)}{q-J-1} 
\end{eqnarray*}
The last equality follows from Lemma~\ref{Lemnu}, and we have Equation~\ref{Eqwsinf} of Theorem~\ref{ThmPops}.

For gaps the length $J=1$ and this becomes
$$ w_{g,1}(\infty) \; = \; \prod_{2 < q, \, q | g} \frac{q-1}{q-2}.$$

\vspace{0.125in}

We have now established all of the claims in Theorem~\ref{ThmPops}.  As corollaries we have the following analogues within the cycles of gaps $\pgap(p^\#)$ among
$p$-rough numbers, for major conjectures about the gaps between prime numbers:
\begin{itemize}
\item {\em k-tuple result} \cite{HL}:  Every admissible constellation $s$ arises in $\pgap(p^\#)$ as $p$ grows.  If the length of $s$ is $J$, then the population of
$s$ grows as 
$$n_{s,J}(p^\#) \; = \; \Theta \left(\prod_{J+1 < q \le p}(q-J-1)\right).$$  
The asymptotic relative population of $s$, compared to other admissible constellations of length $J$,
depends only on the factors $(q-\nu_s(q))$ over the odd prime divisors $Q(s)$ of spans within $s$.
\item {\em Polignac result} \cite{Pol}:  For every even number $2n$, the gap $g=2n$ arises in $\pgap(p^\#)$ as $p$ grows.  The population of $g$
grows as $\Theta(\prod_{2 < q \le p}(q-2))$, and its asymptotic relative population $w_{g,1}(\infty)$ compared to other gaps depends only on the prime factors of $g$.
\item {\em Hardy-Littlewood estimates} \cite{HL}: The asymptotic relative populations $w_{s,J}(\infty)$ for admissible constellations, from Equation~\ref{Eqwsinf}, 
are consistent with
the estimates in Hardy \& Littlewood's Theorem~X1 and Conjecture~B.  A detail in these two conjectures by Hardy \& Littlewood is that Conjecture~B (on
\cite{HL} page 42) is about prime gaps, and Theorem~X1 (on \cite{HL} page 61) is for prime differences.  So from our viewpoint the estimates under Theorem~X1
would include general instances of the constellation.
\end{itemize}

\section{Surviving Eratosthenes sieve}
The results above apply to the cycles of gaps $\pgap(p^\#)$ between the $p$-rough numbers, those
numbers left after Eratosthenes sieve has advanced through the prime $p$.  Which of these gaps survive 
further stages of the sieve to be confirmed as gaps between primes?

We know that the smallest remaining composite number after Eratosthenes sieve has advanced through $p_k$ is
$p_{k+1}^2$.  All of the gaps in $\pgap(p_k^\#)$ from $p_{k+1}$ and short of $p_{k+1}^2$ survive as gaps between primes.  

To make this precise we associate the gaps $g_j$ in $\pgap(p_k^\#)$ with the corresponding $p_k$-rough numbers.
$$ \gamma_0=1, \gap \gamma_1= 1+g_1 = p_{k+1}, \gap \gamma_2 = 1+\gamma_1+\gamma_2 = p_{k+2}, \gap \ldots, \gap \gamma_j = 1 + \sum_{i=1}^j g_i $$

The first gap covers the interval from $1$ to $p_{k+1}$.  Let $g_j$ be the gap just before the $p_k$-rough number $\gamma_j = p_{k+1}^2$.
Then all of the gaps $g_2, \ldots, g_{j-1}$ survive the sieve to be confirmed as gaps among primes.
Under the recursion $\pgap(p_k^\#)\fto \pgap(p_{k+1}^\#)$, the fusion $g_j + g_{j+1}$ will occur.

We say that a gap $g_i$ is in the interval $(a,b)$ iff the corresponding $p_k$-rough number satisfies $a < \gamma_i < b$.
So in $\pgap(p_k^\#)$ all of the gaps in the interval $(p_{k+1}, \; p_{k+1}^2)$ survive further stages of the sieve to be confirmed as gaps among primes.

On the other hand, all of the gaps in the interval $(p_{k+1}, \; p_k^2)$ were confirmed as gaps between primes by previous stages of the sieve.  
So the interval $(p_k^2, \, p_{k+1}^2)$
is the unique interval of gaps in the cycle $\pgap(p_k^\#)$ first confirmed as gaps between primes by the stage of the sieve for $p_k$.  

We call this interval the {\em interval of survival} 
$$\Delta H(p_k) \; = \; (p_k^2, \, p_{k+1}^2).$$
We note that $\Delta H(p_k)$ includes the gap across $p_k^2$ and not the gap across $p_{k+1}^2$.
Both endpoints are composite numbers, and we have to assign the gaps across them in a consistent way.
By definition the gap crossing $p_k^2$ is included in $\Delta H(p_k)$, while the gap crossing $p_{k+1}^2$ is not.

\vspace{0.125in}

We make the following conjecture to support initial estimates of survival.

\begin{conjecture}\label{ConjUni}
The populations of gaps and constellations in the interval of survival $\Delta H(p_k)$ are an approximately uniform sampling of the gaps in $\pgap(p_k^\#)$.
\end{conjecture}

\vspace{0.1in}

Under Conjecture~\ref{ConjUni} we would expect the populations of constellations within $\Delta H(p_k)$ to agree principally
with the relative populations $w_{s,J}(p_k^\#)$.   With ${g=p_{k+1}-p_k}$, we make a simple initial estimate. 
\begin{eqnarray}
E_s(p_k) & = & \frac{p_{k+1}^2-p_k^2}{p_k^\#} \cdot n_{s,J}(p_k^\#)  \gap = \gap  \frac{g \cdot (2p_k+g)}{p_k^\#} \cdot n_{s,J}(p_k^\#) \nonumber \\
   & = & \frac{g \cdot (2p_k+g)}{p_k^\#} \cdot w_{s,J}(p_k^\#) \cdot \prod_{J+1 < q \le p_k} (q-J-1) \nonumber \\
   & = & \frac{g \cdot (2p_k+g)}{\prod_{2 \le q \le J+1} q} \cdot w_{s,J}(p_k^\#) \cdot \prod_{J+1 < q \le p_k} \frac{q-J-1}{q} \label{EqDHest}
\end{eqnarray}

The intuition behind this hypothesis is that under the recursion 
${\pgap(p_k^\#)\rightarrow \pgap(p_{k+1}^\#)},$
the $p_{k+1}$ images of each instance created in step R2 are in fact uniformly distributed at distances $p_k^\#$.
The fusions in step R3 do remove $J+1$ of these images, but for ${p >> J}$ these fusions will be a secondary effect.
The fusions also occur symmetrically across $\pgap(\pml{p})$.  The cycle $\pgap(p^\#)$ is symmetric:
in $\pgap(p^\#)$, we have $ g_i \; = \; g_{\phi(p^\#)-i}$.

There is no part of the discrete dynamic system that is biased for or against any particular admissible constellation.
All admissible constellations of the same length grow eventually at the same rate, under the same dynamics.  
The relative populations already
account for the evolutions of different admissible constellations of the same length, across stages of the sieve.

\begin{figure}[hbt]
\centering
\includegraphics[width=6in]{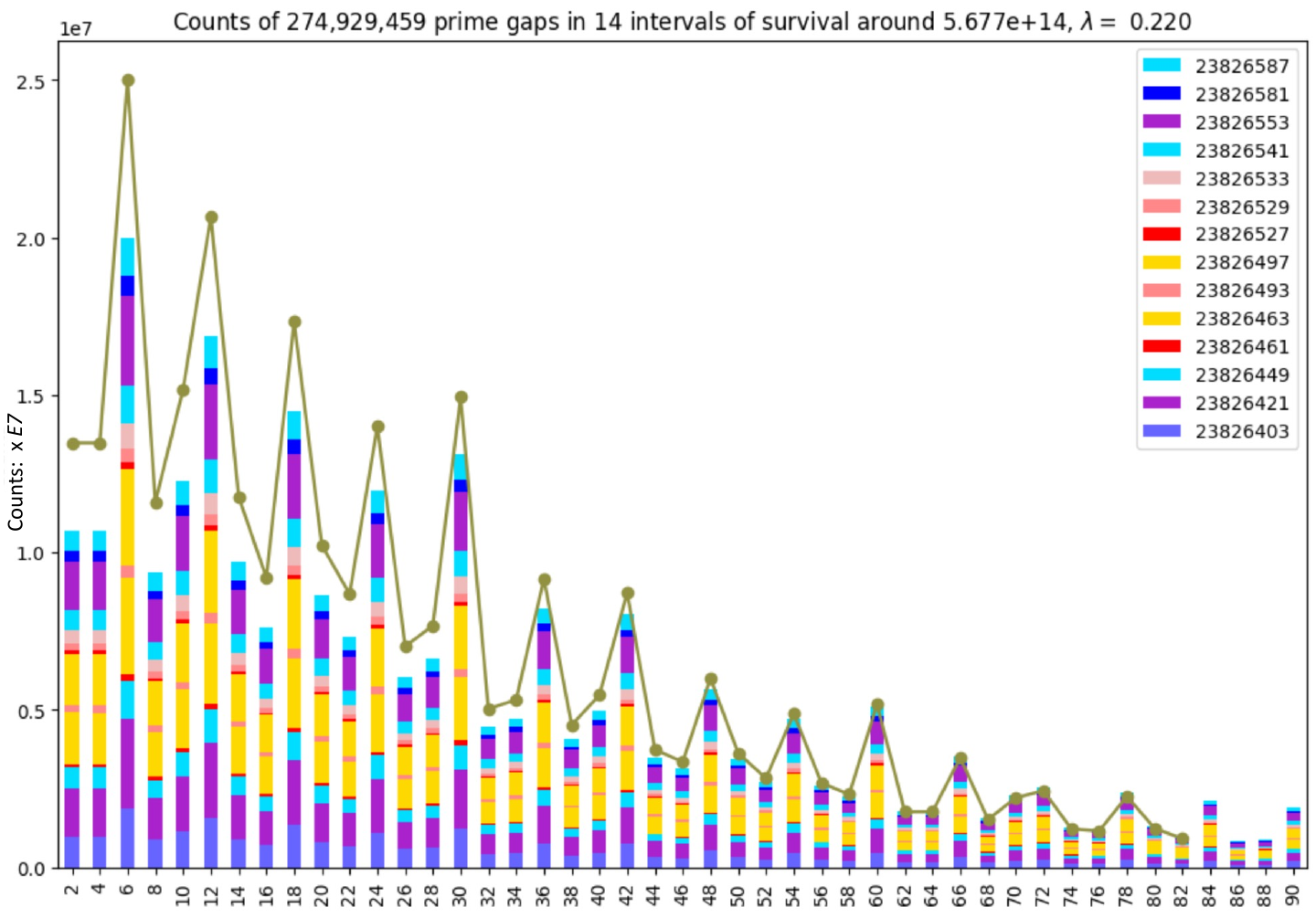}
\caption{\label{DH220Fig} Counts of the gaps between primes in $14$ intervals of survival $\Delta H(p_k)$, for $p_k \approx 2.38\,E7$. 
The individual samples are color-coded by the size of the gap $g=p_{k+1}-p_k$.  The initial estimate in Equation~\ref{EqDHest}
is shown as the broken line. }
\end{figure}

Consequently, as $p$ gets large we would expect the counts of constellations in the interval of survival $\Delta H(p)$
to reflect the relative populations $w_{s,J}(p^\#)$ at this stage of the sieve.
Figure~\ref{DH220Fig} compares the actual populations of the gaps $2 \le g \le 82$ to $E_g(p_k)$ across the $14$ intervals of
survival from $\Delta H(23826403)$ to $\Delta H(23826587)$.

For samples we count the occurrences of gaps in consecutive intervals of survival $\Delta H(p_k)$.  One example of these samples is shown in Figure~\ref{DH220Fig}.
Other samples for gaps are available on GitHub in the notebook {\ttfamily 13\_DelHdisplay}, and the samples for a few select constellations are generated by
the notebook {\ttfamily 22\_DelHsdisplay}.

In the estimate \ref{EqDHest} the only factor that varies with the specific constellation $s$ of length $J$ is the relative population $w_{s,J}(p_k^\#)$.
Fixing a constellation $s$ or gap $g$, over a range of consecutive primes $p_k$ most factors in the estimate \ref{EqDHest} have little variation, except
for the factor ${g=p_{k+1}-p_k}$.   This proportionality to $g$ is reflected as the dominant fluctuation within each column of Figure~\ref{DH220Fig}.

\subsubsection{The missing factor $2e^{-\gamma}$.}
The material up through Section~\ref{ModelSection} and Theorem~\ref{ThmPops} are established independently of and prior to these studies of survival.
Based on the models of relative populations $w_{s,J}(p_k^\#)$ established earlier, we are here formulating estimates of what populations
of gaps and constellations we might expect in the interval $\Delta H(p_k)$.

Our initial estimate $E_s(p_k)$ exhibits the same type of error that Hardy \& Littlewood discuss regarding Sylvester's formula and one of Brun's estimates
(\cite{HL}, p.32--38).  If we take the aggregate $\sum_g E_g(p_k)$ of our initial estimates and compare this total to the total expected under the prime number theorem,
we are off by a factor of $2e^{-\gamma} \approx 1.123$.  That is, in aggregate our estimates $E_g(p_k)$ will be at best about $12\%$ too high.

In our studies we have at times simply introduced a corrective factor to our estimates.  For the present exposition, we have decided to use the
initial estimate that seems to be intrinsically suggested by the models $w_{g,1}(p_k^\#)$.   At this point, with this choice of initial estimate, 
we leave the aggregate $+12\%$ error as a secondary effect.

\subsection{Quadratic density $\eta$}
We introduce a more stable statistic, the {\em quadratic density} $\eta_s(p_k)$ for the occurrences of the constellation $s$ across 
the interval of survival $\Delta H(p_k)$.
The statistic $\eta_s$ and its initial estimate $\widehat{\eta}_s$ are given by
$$\eta_s(p_k) \; = \; N_s(p_k)/ g \biggap {\rm and} \biggap \widehat{\eta}_s(p_k) \; = \; E_s(p_k) / g$$
where $g = p_{k+1}-p_k$.
The sampled values for $\eta_s(p_k)$ are given by the counts $N_s(p_k)$ of occurrences of $s$ within the interval $\Delta H= (p_k^2, p_{k+1}^2)$, divided by $g$.
The estimated quadratic density $\widehat{\eta}_s(p_k)$ is the expected average population of $s$ across the $g$ quadratic intervals $(n^2, (n+1)^2)$ for 
$p_k \le n < p_{k+1}$.

\begin{figure}[hbt]
\centering
\includegraphics[width=6in]{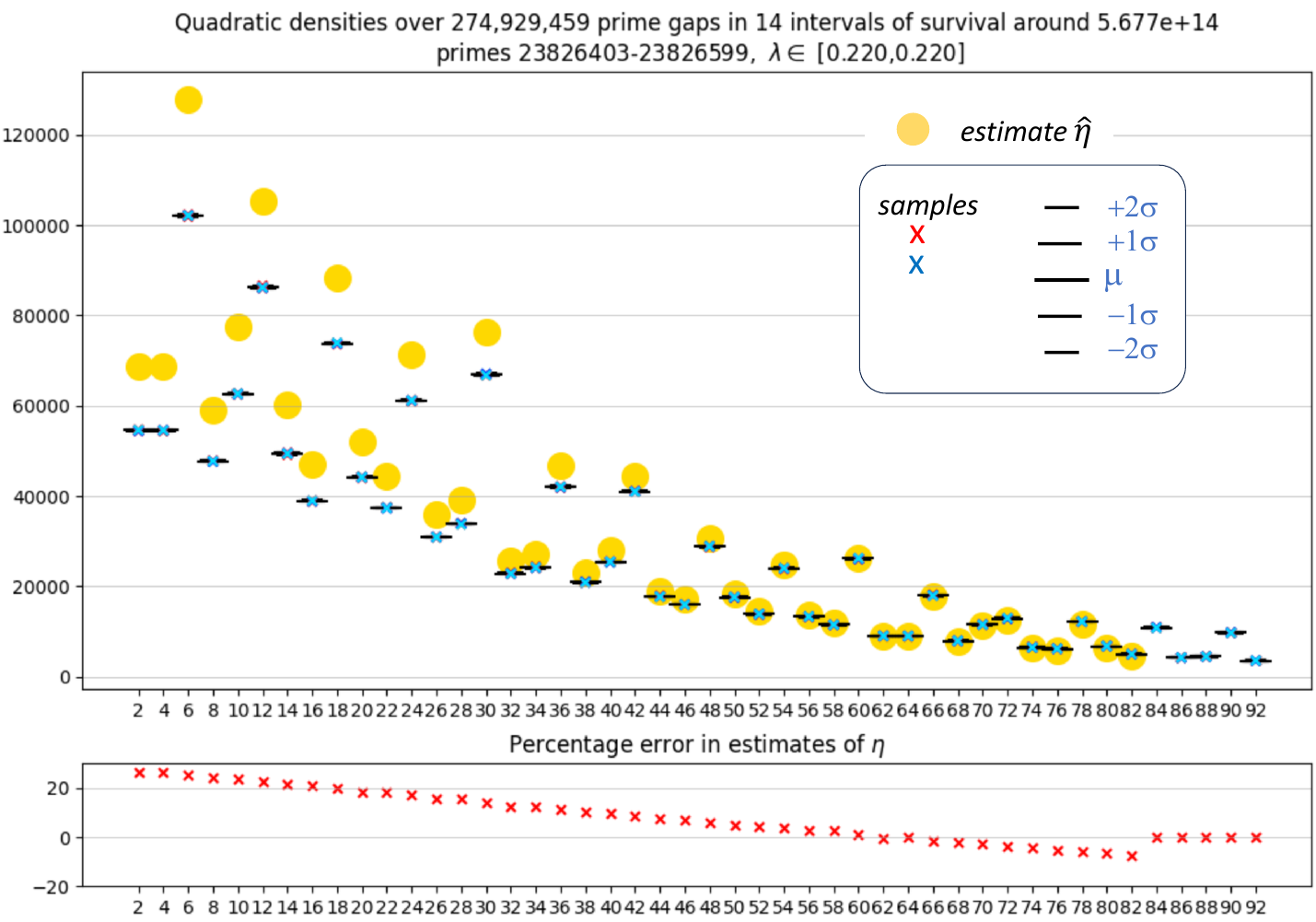}
\caption{\label{eta220Fig} For the same raw counts $N_g(p_k)$ graphed in Figure~\ref{DH220Fig} we here plot the quadratic density samples $\eta_g(p_k)$
and estimates $\widehat{\eta}_g(p_k)$. Within each column the $14$ sampled values are indistinguishable at this scale.  }
\end{figure}

For the raw counts $N_g(p_k)$ shown in Figure~\ref{DH220Fig} we plot the quadratic densities $\eta_g(p_k)$ in Figure~\ref{eta220Fig}.  Here
${23826403 \le p_k \le 23826587}$, and the counts include almost $275$~million prime gaps around $5.677\,E14$.  The estimated quadratic density
$\widehat{\eta_g}(p)$ for each gap is plotted as a gold disk.

\begin{figure}[hbt]
\centering
\includegraphics[width=6in]{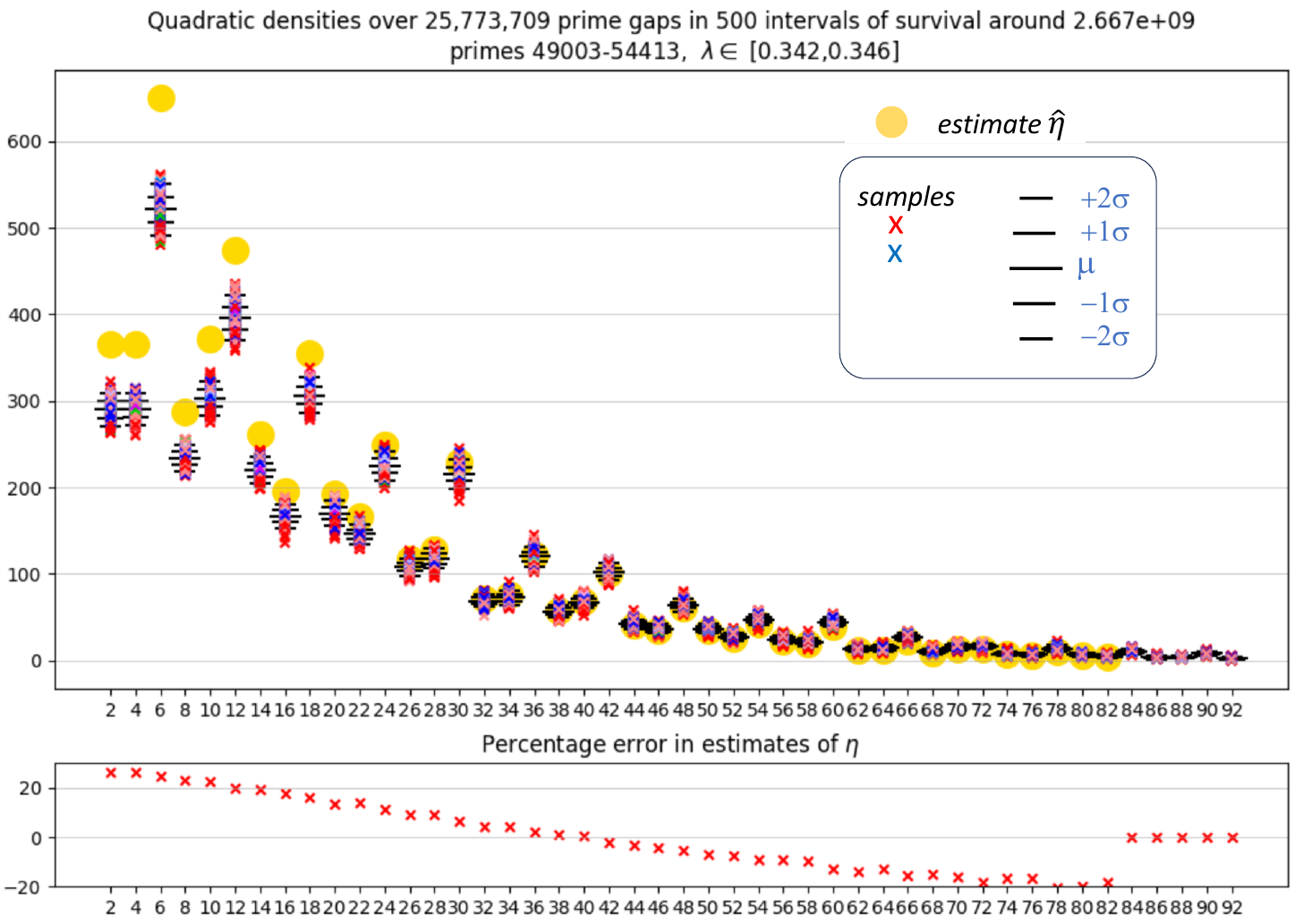}
\caption{\label{eta342Fig} Shown here are the quadratic density samples $\eta_g(p_k)$ and estimates $\widehat{\eta}(p_k)$ for $500$
consecutive intervals of survival $\Delta H(p_k)$, for $49003 \le p_k < 54413$.  For the columns for the smaller gaps, we can see some spread
in the sampled data. }
\end{figure}

The sampled values for the quadratic densities in Figure~\ref{eta220Fig} have low variance.  
We include an earlier sample in Figure~\ref{eta342Fig} across more intervals $\Delta H(p)$
to see the
distributions of the sampled values $\eta_g(p_k)$ within each column.  Figure~\ref{eta342Fig} shows $\eta_g(p_k)$ and the estimate $\widehat{\eta}_g(p_k)$
for the $25$ million prime gaps across $500$ consecutive intervals of survival $\Delta H(p_k)$ for $49003 \le p_k \le 54413$.

Looking at Figures~\ref{eta220Fig}~\&~\ref{eta342Fig}, we observe that our initial estimate $\widehat{\eta}_g$ does generally track the growth of the populations of gaps
in $\Delta H(p)$.  The lower panels in those two figures show the relative error
in the estimate $\widehat{\eta}_g(p)$, and the residual error skews systematically toward the smaller gaps.  

An open question is about secondary effects: why do the sampled values for a certain gap deviate from the initial estimate, which assumes a uniform distribution?
Since the initial estimate and the residual error exhibit a lot of structure, what are the mechanisms that account for these variances?

For gaps we have the initial conditions from $\pgap(37^\#)$, so we have models $w_g(p_k^\#)$ only for gaps $g \le 82$.  
For constellations $s$ for $J \ge 2$ we currently use $\pgap(29^\#)$ for initial conditions.  So we can determine the coefficients $\ell_j$
for the exact models of relative populations $w_{s,J}(p_k^\#)$ for constellations $s$ with span $|s|\le 62$.

The quadratic density $\eta_s(p)$ balances different aspects of the population of an admissible constellation $s$ among primes.  The raw population $n_s(p^\#)$
across the cycles of gaps $\pgap(p^\#)$ grows superexponentially by factors of $(p-J-1)$.  However, $\pgap(p^\#)$ has length $\phi(p^\#)$ and
the overall density of $s$ in $\pgap(p^\#)$ fades to $0$.  The interval of survival ${ \Delta H(p_k) =(p_k^2, p_{k+1}^2)}$ segments the surviving instances of $s$ among
primes, but as evidenced in Figure~\ref{DH220Fig} and Equation~\ref{EqDHest}, the behavior of the counts of surviving instances of $s$ in $\Delta H(p_k)$
is dominated by the fluctuations of the gap ${g=p_{k+1}-p_k}$.

The quadratic density $\eta_s(p_k)$ provides a more stable statistic for the counts and
estimates of $s$ that occur between primes in $\Delta H(p_k)$, and it has the intuitive interpretation that $\widehat{\eta}_s(p_k)$ is the expected number of occurrences
of $s$ between primes in each of the quadratic intervals $(n^2, (n+1)^2)$ across $\Delta H(p_k)$.

\subsection{When does quadratic density increase?}
We look for trends in the large-scale behavior of the estimate $\widehat{\eta}_s$.  In Lemma~\ref{fLemma} below we identify a condition under which
the expected quadratic density $\widehat{\eta}_s$ increases from $\Delta H(p_{k-1})$ to $\Delta H(p_k)$.  For the specific case of gaps, 
for which $J=1$, Theorem~\ref{gapThm} below will show
that the expected quadratic density for any gap $g$ always increases once $g$ arises in the cycle $\pgap(p^\#)$.

\vspace{0.125in}

We want to identify conditions under which $\widehat{\eta}_s(p_k) > \widehat{\eta}_s(p_{k-1})$.

\begin{figure}[hbt]
\centering
\includegraphics[width=2.75in]{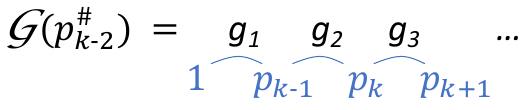}
\caption{\label{g2g3Fig} To articulate the condition under which ${\widehat{\eta}(p_{k-1}) < \widehat{\eta}(p_k)}$, we need the gaps
$g_2$ from $p_{k-1}$ to $p_k$ and $g_3$ from $p_k$ to $p_{k+1}$.}
\end{figure}

To set the stage for this analysis, let's start in the cycle $\pgap(p_{k-2}^\#)$.  The first gap $g_1=p_{k-1}-1$.  The second gap $g_2 = p_k-p_{k-1}$,
and the third gap $g_3=p_{k+1}-p_k$.  See Figure~\ref{g2g3Fig}.

For an admissible constellation $s$ of length $J$, let $p_0$ be large enough that $s$ occurs in the cycle $\pgap(p_0^\#)$ and $|s|< 2p_1$.
For $p_k > p_0$, 
$$ n_{s,J}(p_k^\#) \; = \; (p_k-J-1) \cdot n_{s,J}(p_{k-1}^\#) \, + \, n_{s,J+1}(p_{k-1}^\#),$$
as illustrated in Figure~\ref{DynFig}.   Let $\delta_s$ be the ratio of the number of driving terms for $s$ of length $J+1$
in the cycle $\pgap(p_{k-1}^\#)$ to the number of instances of $s$ itself in that cycle. 
$${\delta_s = n_{s,J+1}(p_{k-1}^\#) \; / \; n_{s,J}(p_{k-1}^\#)}$$ 

We define
\begin{equation} 
f(p_k, J) \; = \; \frac{g_2+g_3}{2+ g_3/p_k} - (J+1). \label{Eqf}
\end{equation}
This function $f$ is linear in $J$ and is independent of the constellation $s$.

We are ready to state the condition under which $\widehat{\eta}_s(p_k)$ is increasing.

\begin{lemma}\label{fLemma}
Suppose the admissible constellation $s$ of length $J$ occurs in the cycle of gaps $\pgap(p_0^\#)$, and $|s| < 2p_1$.
Then for $p_k > p_0$,
\begin{eqnarray}
\widehat{\eta}_s(p_k) & >  & \widehat{\eta}_s(p_{k-1}) \nonumber \\
 &  {\rm iff} &  g_2 + g_3 \; > \; (J+1 - \delta_s) \cdot (2 + \frac{g_3}{p_k}) \label{Eqg2g3} \\ 
  & {\rm iff} &  f(p_k,J)+\delta_s \; > \; 0 \label{Eqfdel}
 \end{eqnarray}
 \end{lemma}

\begin{proof}
\begin{eqnarray*}
\frac{\widehat{\eta}_s(p_k)}{\widehat{\eta}_s(p_{k-1})} & = & \frac{E_s(p_{k})/ g_3}{E_s(p_{k-1})/g_2} \; = 
  \; \frac{n_{s,J}(\pml{p_k})}{n_{s,J}(\pml{p_{k-1}})} \cdot \frac{p_{k-1}^\#}{p_k^\#}
  \cdot \frac{g_3 (p_{k+1}+p_k)/g_3}{g_2 ( p_k+p_{k-1})/g_2} \\
  & = &  \frac{1}{p_k} \cdot \frac{ (p_{k+1}+p_k)}{( p_k+p_{k-1})} \cdot \frac{n_{s,J}(\pml{p_k})}{n_{s,J}(\pml{p_{k-1}})} \\
  & = & \frac{1}{p_k} \cdot \frac{ (p_{k+1}+p_k)}{( p_k+p_{k-1})} \cdot \frac{(p_k-J-1) n_{s,J}(\pml{p_{k-1}})+n_{s,J+1}(\pml{p_{k-1}})}{n_{s,J}(\pml{p_{k-1}})}  \\
  & = & \frac{1}{p_k} \cdot \frac{ (p_{k+1}+p_k)}{( p_k+p_{k-1})} \cdot \left[ (p_k-J-1) + \frac{n_{s,J+1}(\pml{p_{k-1}})}{n_{s,J}(\pml{p_{k-1}})} \right] \\
  & = & \frac{1}{p_k} \cdot \frac{ (p_{k+1}+p_k)}{( p_k+p_{k-1})} \cdot \left[ (p_k-J-1) + \delta_s \right] 
\end{eqnarray*}
Thus $\frac{\widehat{\eta}_s(p_k)}{\widehat{\eta}_s(p_{k-1})} > 1$ iff
\begin{align}
\frac{1}{p_k} \cdot \frac{ (p_{k+1}+p_k)}{( p_k+p_{k-1})} \cdot \left[ (p_k-J-1) + \delta_s \right] & \; >  1 \nonumber \\
\frac{ (p_{k+1}+p_k)}{( p_k+p_{k-1})} - 1 & \; >  \frac{1}{p_k} \cdot \frac{ (p_{k+1}+p_k)}{( p_k+p_{k-1})} \cdot \left[ J+1 - \delta_s \right]  \nonumber \\
 (p_{k+1} -p_{k-1}) & \; >  \frac{1}{p_k} \cdot (p_{k+1}+p_k) \cdot \left[ J+1 - \delta_s \right]  \nonumber \\
 g_2 + g_3 & \; >  \frac{1}{p_k} \cdot (2 p_k + g_3) \cdot \left[ J+1 - \delta_s \right]  \nonumber \\
 g_2 + g_3 & \; >  (J+1 - \delta_s) \cdot (2 + \frac{g_3}{p_k})  \tag{\ref{Eqg2g3}}
\end{align}
This is the desired middle inequality \ref{Eqg2g3} in the lemma.
Continuing, we have
\begin{eqnarray*}
g_2 + g_3 & > & (J+1 - \delta_s) \cdot (2 + \frac{g_3}{p_k}) \\
\frac{g_2 + g_3}{2 + g_3/p_k} & > & J + 1 - \delta_s \\
f(p_k, J) + \delta_s & > & 0
\end{eqnarray*}
\end{proof}

$f$ captures the part of the constraint~(\ref{Eqg2g3}) that does not vary with the specific constellation $s$, and (\ref{Eqfdel}) isolates the contribution of $s$ as 
a simple offset $+\delta_s$.  
The ratio $\delta_s\ge 0$ with equality iff $n_{s,J+1}(p_{k-1}^\#) = 0$.  Thus whenever $f (p_k,J) > 0$ the estimated quadratic density increases for {\em all} admissible 
constellations $s$ in $\pgap(\pml{p_k})$ of length $J$ and $|s| < 2p_k$.

\begin{figure}[hbt]
\centering
\includegraphics[width=5.75in]{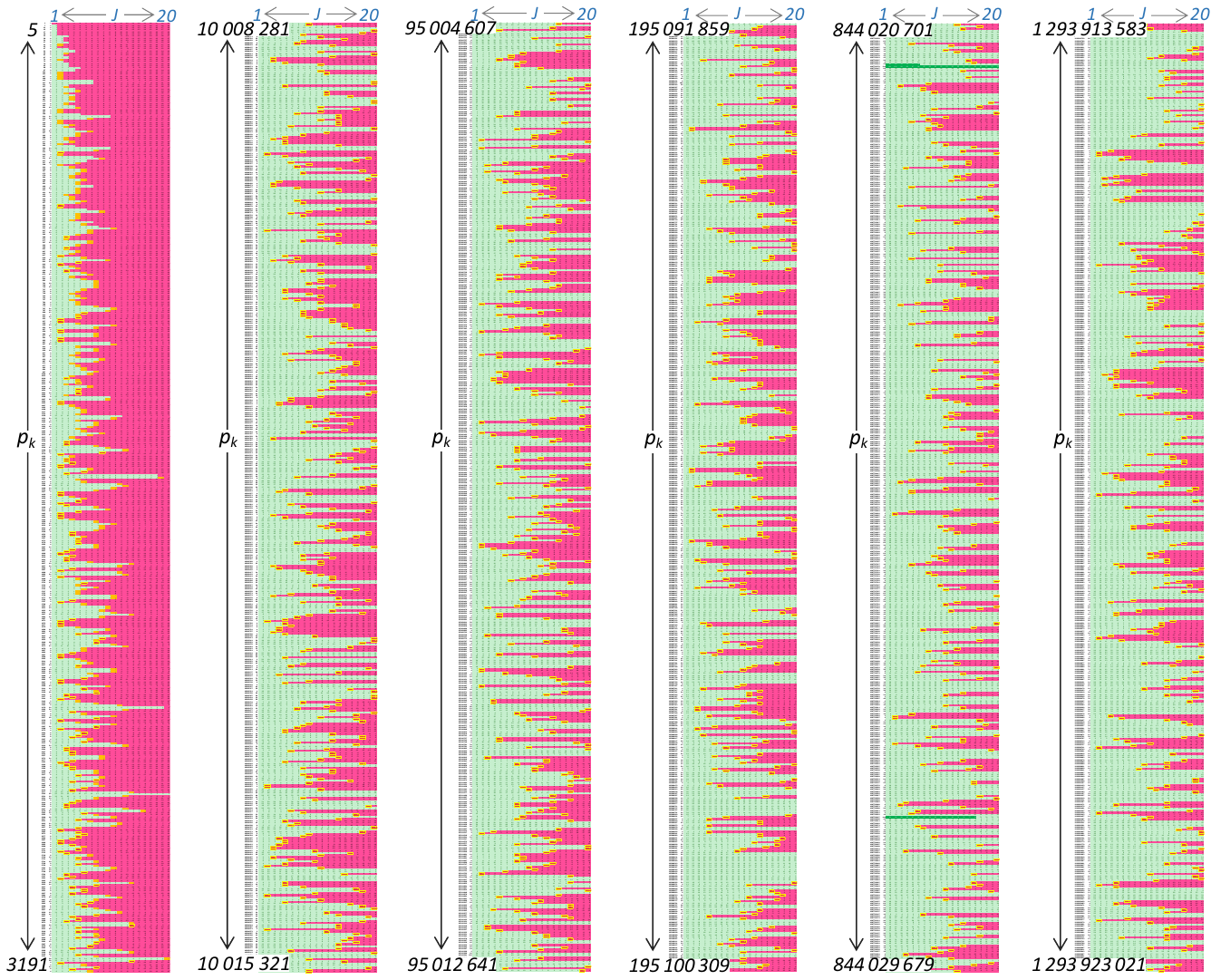}
\caption{\label{fpjFig} Six samples of $f(p_k,J)$ over consecutive primes and lengths ${1\le J \le 20}$.  
When ${f(p,J)>0}$ (green) the expected quadratic density ${\widehat{\eta}_s(p,J)}$
is increasing.  When ${f(p,J)<0}$ (red), $\widehat{\eta}_s(p,J)$ is decreasing unless the ratio $\delta_s(p)$ makes up the deficit.
Across each row the value of $f$ decreases linearly with $J$.  
The variations in $f(p,J)$ across the rows reflects the value of ${g_2+g_3 = p_{k+1}-p_{k-1}}$.  }
\end{figure}

Figure~\ref{fpjFig} shows $f(p,J)$ over six samples of consecutive primes within the first $65$~million primes.  
For these samples we show ${1 \le J \le 20}$.
Across each row, $f(p,J)$ is linear in $J$.  The pairs $(p,J)$ with $f(p,J) > 0$ are colored green, and the pairs with $f(p,J) < 0$ are colored red.
When $g_2+g_3$ is small, the sign change in $f$ moves toward smaller $J$'s, and when $g_2+g_3$ is larger, the sign change in $f$ moves to
larger values of $J$.

When $f(p_k,J) <0$, the constraint is not met -- unless the addition of $\delta_s$ is sufficient to change the sign for the constellation $s$.  
Early in the evolution of the systems
in Figure~\ref{DynFig} the ratio $\delta_s$ can be quite large.  As primes get large, for any 
admissible constellation $s$ the ratio $\delta_s \fto 0$.  So for large primes the condition $f(p_k,J) \ge 0$ provides a surrogate for the 
constraint~(\ref{Eqfdel}) as $\delta_s \fto 0$.

Note that in Lemma~\ref{fLemma} the constraint~(\ref{Eqg2g3}) holds if and only if $\widehat{\eta}_s$ is increasing at $p_k$.  Expressed as it is in
Equation~(\ref{Eqg2g3}) the constraint ties the behavior of $\widehat{\eta}_s$ simply and explicitly to the succession of gaps $g_2, g_3$, to the length
$J$, qualified surprisingly by the ratio $\delta_s$, and to the ratio $\frac{g_3}{p_k}$.

The lefthand side of the constraint~(\ref{Eqg2g3})  is the span $g_2+g_3$ of the constellation $g_2 g_3$ from $p_{k-1}$ to $p_{k+1}$.  See Figure~\ref{g2g3Fig}.
This is the only occurrence of the gap $g_2$ in the constraint.  If this span $g_2+g_3$ is above the bound given by the righthand side of the constraint,
 then $\widehat{\eta}_{s}$ increases at $p_k$, and otherwise $\widehat{\eta}_{s}$ decreases at $p_k$.

The smallest admissible constellations of length $2$ are $s=2\, 4$ and $s=4\, 2$, each of which has span $|s|=6$, so ${g_2+g_3 \ge 6}$.

\begin{theorem}\label{gapThm}
For every gap $g$, once $g$ occurs in $\pgap(\pml{p_0})$ for any $p_0\ge 3$ with $g < 2p_1$, 
the expected quadratic density $\widehat{\eta}_g$ increases for all primes $p_k > p_1$.
\end{theorem}

\begin{proof}
For a gap $g$, the length $J=1$.  The parameter $\delta_g \ge 0$, and ${0 < \frac{g_3}{p_k} < 1}$. So 
$$ g_2 + g_3 \; \ge \; 6 \;  > \; (2-\delta_g)\cdot (2 + \frac{g_3}{p_k}), $$
and the inequality in (\ref{Eqg2g3}) always holds.  Thus $\widehat{\eta}_g$ is increasing.
\end{proof}

Under the dynamics of Eratosthenes sieve, the estimated quadratic density of every gap increases.
For any gap $g$ and large enough primes $p_k$, the expected population of the gap $g$ in quadratic intervals $(n^2, (n+1)^2)$ grows from
one interval for survival $\Delta H(p_k)$ to the next.

For $J \ge 2$, results along the lines of Theorem~\ref{gapThm} come with qualifications.  
For example, for length $J = 2$, the constraint~(\ref{Eqg2g3})
will hold for large $p_k$ when $g_2 + g_3 \ge 8$, or equivalently except when 
the constellation $g_2g_3 = 2\, 4$ or $4\, 2$.

On the other hand, $g_2+g_3$ is the most volatile part of the constraint (\ref{Eqg2g3}).  Its minimum value is $6$, but the value $g_2+g_3$ 
fluctuates wildly.  Look at the samples for $f(p,J)$ in Figure~\ref{fpjFig}.  The mean gap size between primes is growing, but successions of small gaps continue
to occur.

Results along the line of Theorem~\ref{gapThm} for longer constellations will also be sensitive to the value of ${\delta_s= \frac{n_{s,J+1}(\pml{p_k})}{n_{s,J}(\pml{p_k})}}$,
 which depends both on the
constellation $s$ and the prime $p_k$.  

\section{Conclusion}
We approach Eratosthenes sieve as a discrete dynamic system.  The objects of the system are the cycles of gaps $\pgap(p^\#)$
at each stage of the sieve.  There is a 3-step recursion that produces the next cycle of gaps from the current one.
$$ \pgap(p_k^\#) \; \fto \; \pgap(p_{k+1}^\#).$$
A lot of structure in $\pgap(p^\#)$ is preserved under this recursion.

Taking initial conditions from the cycle $\pgap(p_0)$, we have extracted exact models for the relative populations
of all admissible constellations $s$ for which the span $|s| \le p_1$.  These models allow us to compare relative populations among
admissible constellations far beyond the computational horizon.  For computed samples of the populations of constellations, we can compare 
the samples to the values expected from the models.

In order to study the correlation between the populations of admissible constellations in $\pgap(p_k^\#)$ and the constellations of gaps among primes,
we define the interval of survival $\Delta H(p_k) = (p_k^2, p_{k+1}^2)$ and introduce a statistic $\eta_s(p_k)$, the quadratic density of $s$ in 
$\Delta H(p_k)$.  The quadratic density $\eta_s(p_k)$ is the average number of occurrences of $s$ in each quadratic interval $(n^2, (n+1)^2)$
across $\Delta H(p_k)$.  We compare sampled values of $\eta_s(p_k)$ to estimates $\widehat{\eta}_s(p_k)$ derived from the relative population
$w_{s,J}(p_k^\#)$.

Ongoing avenues of research on this project include the following:
\begin{itemize}
\setlength{\itemsep}{0.02in}
\item Develop the coefficients for the models $w_s(p^\#)$ for more gaps and constellations.  Tomasz Hojka \cite{THj} has recently calculated the coefficients
$\ell_j$ for the models $w_{g,1}(p^\#)$ for the gaps $84 \le g \le 104$.  His code and manuscripts are available at 
\mbox{\url{https://github.com/hojka1987t-dot/sieve-dynamics-beyond-2p}}.
\item Take more samples $\eta_s(p_k)$ for gaps and constellations for intervals of survival $\Delta H(p_k)$ for primes $p_k > 4\,E7$.
\item Identify a secondary correction to the estimated quadratic density $\widehat{\eta}_s(p_k)$, at least for small gaps. 
\end{itemize}

The code and data supporting the analysis above can be found in Jupyter notebooks at \mbox{\url{https://github.com/fbholt/Primegaps-v2}}.

\noindent{\em Acknowledgements.}
The author is grateful to the anonymous reviewers for providing several helpful suggestions, to make the exposition clearer and more precise.


\begin{thebibliography}{HR15}

\normalsize
\baselineskip=17pt

\bibitem[HL23]{HL} G.H. Hardy and J.E. Littlewood,
\emph{Some problems of 'Partitio numerorum' {III}: on the expression of a number as a sum of primes},
Acta Math. 44 (1923), 1-70, in G.H.~Hardy Collected Papers, Clarendon Press (1966).

\bibitem[Hoj26]{THj} T. Hojka,
\emph{Exact dynamics of {E}ratosthenes sieve beyond $2p$: coincident fusions, births of gaps and missing differences},
preprint (2026).

\bibitem[HR15]{FBHSFU} F.B. Holt and H. Rudd,
\emph{Combinatorics of the gaps between primes},
Connections in Discrete Mathematics - Simon Fraser U. (2015),
arXiv:1510.00743.

\bibitem[Ho22]{FBHPatterns} F.B. Holt,
Patterns among the Primes, KDP (2022).

\bibitem[Ho25]{FBHktuple} F.B. Holt,
\emph{Eratosthenes sieve supports the $k$-tuple conjecture},
preprint (2025), arXiv:2502.20470.

\bibitem[Ho24]{FBHbias} F.B. Holt,
\emph{Expected biases in the distribution of consecutive primes},
preprint (2024), arXiv:2405.03540.

\bibitem[Ho23]{FBH2p} F.B. Holt,
\emph{Models for gaps $g=2p_1$},
preprint (2023), arXiv:2309.16833v2.

\bibitem[OS16]{OS} R. Lemke Oliver and K. Soundararajan,
\emph{Unexpected biases in the distribution of consecutive primes},
Proc. Natl. Acad. Sci USA 113 (2016), 4446-54.

\bibitem[dP49]{Pol} A. de Polignac,
\emph{Six propositions arithmologiques of d\'{e}duites du crible d'Eratosth\`{e}ne},
Nouvelles annales de math\'{e}matiques $1^{re}$ s\'{e}rie 8 (1849), 423-29.

\bibitem[Pr82]{Wheel} P. Pritchard,
\emph{Explaining the wheel sieve},
Acta Informatica 17 (1982), 477-485.

\end{thebibliography}
\end{document}